\documentclass[a4wide,11pt,onetabnum,reqno]{article} 
\usepackage[utf8]{inputenc}
\usepackage{fullpage} 
\usepackage{authblk}
\usepackage{amsmath, amssymb, amsfonts, amsbsy, amsthm}
\usepackage{mathrsfs}
\usepackage{bm}
\usepackage{graphicx}
\usepackage{graphics}
\usepackage{epsfig}
\usepackage{titlesec}
\usepackage{enumerate}
\usepackage[table]{xcolor}
\usepackage{pst-all}
\usepackage{calligra}
\usepackage{tikz}
\usetikzlibrary{patterns,calc,decorations.markings}
\usepackage{caption}
\usepackage{algorithm}
\usepackage{multirow}
\usepackage{array}
\usepackage{accents}
\usepackage{blindtext}
\usepackage[nottoc]{tocbibind} 
\usepackage{subcaption}

\usepackage{hyperref} 
\usepackage{cleveref}

\newcommand{\Hrule}[3][.]{%
  \par\addvspace{#2}%
  \begingroup\color{#1}%
  \hrule
  \endgroup
  \addvspace{#3}%
}

\newcolumntype{H}{>{\setbox0=\hbox\bgroup}c<{\egroup}@{}}

\DeclareMathAlphabet{\mathpzc}{OT1}{pzc}{m}{it}
\DeclareMathAlphabet{\mathcalligra}{T1}{calligra}{m}{n}

\newtheorem{theorem}{\bf Theorem}[section]
\newtheorem{proposition}[theorem]{\bf Proposition}
\newtheorem{definition}{\bf Definition}[section]

\newtheorem{exam}[theorem]{\bf Example}
\newtheorem{remark}[theorem]{\bf Remark}
\newtheorem{lemma}[theorem]{\bf Lemma}

\newcommand{\be}{\begin{equation}}
\newcommand{\ee}{\end{equation}}
\newcommand{\beno}{\begin{equation*}}
\newcommand{\eeno}{\end{equation*}}
\newcommand{\ba}{\begin{align}}
\newcommand{\ea}{\end{align}}
\newcommand{\bano}{\begin{align*}}
\newcommand{\eano}{\end{align*}}
\newcommand{\bea}{\begin{eqnarray}}
\newcommand{\eea}{\end{eqnarray}}
\newcommand{\beano}{\begin{eqnarray*}}
\newcommand{\eeano}{\end{eqnarray*}}

\numberwithin{equation}{section}

\newcommand{\oldfactorial}[1]{\mathpalette\oldfactorialaux{#1}}
\newcommand{\oldfactorialaux}[2]{%
  {#1\mkern1mu\oalign{\vrule\,$#1#2\mathstrut$\,\cr\noalign{\hrule}}}}

\def \div {\mathrm{div}}

\def  \argmin {\emph{argmin~}}

\def\sjump#1{[\hskip -1.5pt[#1]\hskip -1.5pt]}

\def \meas {\textit{ meas}}

\def \esssup{\textit{esssup}}

\begin{document}

\title{Adaptive finite element convergence analysis of \textsf{AT1} phase-field model for quasi-static fracture in strain-limiting solids\footnote{Dedicated to the cherished memory of Professor K. R. Rajagopal (Professor of Mechanical Engineering, Texas A\&M University, College Station, Texas-USA). His foundational work and invaluable mentorship left an indelible mark on this research and its contributing authors.}}

\author[1]{Ram Manohar}
\author[1,*]{S. M. Mallikarjunaiah}

\affil[1]{Department of Mathematics \& Statistics, Texas A\&M University-Corpus Christi, Texas- 78412, USA}
\affil[*]{Corresponding author}
\affil[ ]{\textit{E-mail addresses:} \texttt{ram.manohar@tamucc.edu} (Ram Manohar), \texttt{M.Muddamallappa@tamucc.edu} (S.M. Mallikarjunaiah)}

\date{} 
\maketitle

\begin{abstract}
This research rigorously investigates the convergence of adaptive finite element methods for regularized variational models of quasi-static brittle fracture in elastic solids. We specifically examine a novel Ambrosio-Tortorelli (AT1) phase-field model within the framework of elasticity theories, particularly for material models characterized by an algebraically nonlinear stress-strain relationship. Two distinct and novel adaptive mesh refinement algorithms, underpinned by robust local error indicators, were introduced to efficiently solve the underlying nonlinear energy minimization problem. A detailed convergence analysis was conducted on the sequences of minimizers produced by these strategies.

Our findings rigorously demonstrate that the minimizer sequences from the first adaptive algorithm achieve convergence to a predefined tolerance. Crucially, the second algorithm is proven to generate inherently convergent sequences, thereby eliminating the need for an explicit stopping criterion. The practical effectiveness of this proposed adaptive framework is thoroughly validated through extensive numerical simulations. A case study involving an edge crack in an elastic body, governed by an algebraically nonlinear strain-limiting relationship and subjected to anti-plane shear-type loading, is presented. Critical comparisons of the energy components—bulk, surface, and total—showcase the superior performance of both adaptive algorithms.
\end{abstract}
\vspace{.1in}

	
\vspace{.1in}

\section{Introduction to brittle fracture models}

A foundational mathematical framework for modeling brittle fracture in materials is provided by the \textit{Francfort-Marigo model} \cite{Francfort1998}. The total energy of the system within this model is comprised of three fundamental components: \textit{elastic energy}, which represents the energy stored due to material deformation; \textit{bulk energy}, about the intrinsic material properties and internal structure; and \textit{surface energy}, associated with the formation of new fracture surfaces. While a material's resistance to deformation is quantified by its elastic energy, the energy required for the creation of new fracture surfaces is captured by the surface energy. These components collectively define the system's energy landscape.

\subsection{Governing equations for strain-limiting elasticity}

Consider a smooth, open, connected, and bounded domain $\mathcal{D}:=\mathcal{D}(t) \subset \mathbb{R}^{d},\, d=2,\,3$, possessing a given boundary $\partial \mathcal{D}$, with $\Gamma(t)$ representing a crack set. This domain is considered to be occupied by an elastic body whose mechanical response is governed by a specialized nonlinear constitutive relation between the Cauchy stress tensor $\boldsymbol{T} \colon \mathcal{D} \to \mathbb{R}^{d \times d}_{sym}$ and the linearized elasticity tensor $\boldsymbol{\epsilon} \colon \mathcal{D} \to \mathbb{R}^{d \times d}_{sym}$. Within the framework of \textit{strain-limiting theories of elasticity} \cite{rajagopal2007elasticity}, this constitutive class of relations is expressed in the form:
\begin{equation}\label{model1}
\mathcal{F}(\boldsymbol{T}, \; \boldsymbol{B}) = \boldsymbol{0},
\end{equation}
where $\mathcal{F} \colon \mathbb{R}^{d \times d}_{sym} \times \mathbb{R}^{d \times d}_{sym} \to \mathbb{R}^{d \times d}_{sym}$ denotes a nonlinear, tensor-valued function, while $\boldsymbol{B}$ represents the left Cauchy-Green strain tensor. Through linearization, under the assumption of small displacement gradients, nonlinear constitutive relations such as the following can be derived \cite{mallikarjunaiah2015direct,itou2017nonlinear,itou2018states,MalliPhD2015,rajagopal2014nonlinear,rodriguez2021stretch}:
\begin{equation}\label{model2}
\boldsymbol{\epsilon} :=\mathcal{F}(\boldsymbol{T}).
\end{equation}
The following system of equations governs the material behavior \cite{manohar2024hp,rajagopal2011modeling,HCY_SMM_MMS2022}:
\begin{align}\label{model_1}
\begin{cases}
\div \, \boldsymbol{T} = \boldsymbol{f}, \quad \boldsymbol{T}=\boldsymbol{T}^T & \mbox{in} \quad \mathcal{D}, \\
\boldsymbol{\epsilon} :=\mathcal{F}(\boldsymbol{T}) & \mbox{in} \quad \mathcal{D}.
\end{cases}
\end{align}
In Equation \eqref{model_1}, the first part of the first equation represents the \textit{balance of linear momentum}, while the second part implies the \textit{symmetry of the stress tensor} in the absence of body couples. For the specific case of anti-plane shear loading, the displacement vector is characterized by only one nonzero component (i.e., in the $z$-direction), which results in planar stress and strain tensors with only a single nonzero component. If $u(x,y)$ is the \textit{Airy's stress function}, satisfying the relations:
\begin{equation}\label{Tcompo}
\boldsymbol{T}_{13} = \dfrac{\partial u}{\partial y}, \quad \boldsymbol{T}_{23} = - \dfrac{\partial u}{\partial x},
\end{equation}
then the equilibrium equation is automatically satisfied. Consequently, the compatibility condition is reduced to:
\begin{equation}\label{sc}
\dfrac{ \partial \boldsymbol{\epsilon}_{13} }{\partial y} - \dfrac{ \partial \boldsymbol{\epsilon}_{23} }{\partial x} =0.
\end{equation}
The current study is specifically devoted to analyzing a subclass of the aforementioned general models. In this work, a material model is considered in which the linearized strain is expressed as a nonlinear function of stress, as given by:

\begin{equation}\label{model_2}
\boldsymbol{\epsilon} := \dfrac{\boldsymbol{T}}{2\mu \left( 1 + \beta^{\alpha} \, |\boldsymbol{T}|^{\alpha} \right)^{1/\alpha}}.
\end{equation}
Here, $\beta \geq 0$ and $\alpha >0$ are defined as modeling parameters, and $\mu>0$ represents the shear modulus. A similar constitutive class has been investigated in previous studies \cite{gou2023MMS,gou2023computational,kulvait2013,mallikarjunaiah2015direct,vasilyeva2024generalized,HCY_SMM_MMS2022,yoon2022finite,yoon2021quasi}. Substituting equations \eqref{Tcompo} into \eqref{model_2}, and subsequently utilizing the strain components within the compatibility condition \eqref{sc}, yields a second-order quasi-linear elliptic partial differential equation.

\subsection{Evolution of fracture mechanics and phase-field models}
The theoretical aspects of quasi-static and brittle fracture mechanics were first considered by Francfort and Marigo \cite{Francfort1998}, who also successfully conducted experiments on a two-dimensional example involving a domain with a straight initial crack. By rephrasing Griffith's concept \cite{Griffith1921} of balancing energy release rate with a fictitious surface energy as an energy minimization problem; Francfort and Marigo developed a model that circumvented the standard constraints of classical fracture mechanics, such as a predetermined and piecewise smooth crack path. The Francfort-Marigo model approach is predicated on the minimization of a highly irregular energy functional, also known as the Mumford-Shah function in the context of image segmentation \cite{Mumford1989}. To render this problem amenable to numerical simulation, numerous regularization techniques have been introduced in the literature \cite{BourdinChambolle2000, Chambolle1999, Negri2001}. These techniques typically construct approximation functionals whose minimizers are demonstrated to converge to those of the original functional through the notion of $\Gamma$-convergence \cite{Braides2002}.

Furthermore, the Ambrosio-Tortorelli approximation \cite{Ambrosiott1990, Ambrosiottapprox1992} possesses particularly attractive properties among regularization methods, notably that its minimization can be reduced to the solution of elliptic boundary value problems, which are straightforward to discretize, for instance, using a finite element approach. This strategy has been effectively employed by Bourdin, Francfort, and Marigo \cite{BourdinChambolle2000} and Bourdin \cite{Bourdin2007, Bourdinvariational2007} for simulating problems that typically fall outside the scope of conventional techniques. The Ambrosio-Tortorelli approximation can be conceptualized as a phase-field model for the fracture set. To accurately resolve the phase-field variable, the mesh in the vicinity of the fracture must be significantly finer than that required for determining the elastic deformation alone. Consequently, the utilization of an adaptively refined mesh becomes a natural consideration, particularly as the fracture path is not known \textit{a priori}.

\subsection{Our Contribution.}
The entire paper is arranged as follows: In \textbf{Section 2}, function spaces, basic definitions, and the nonlinear setting of the Francfort-Marigo model of brittle fracture are introduced. Using these fundamental arguments and explanations, the Fr\'echet differentiability of the minimization functional is demonstrated. Furthermore, it is proven that every phase-field minimizer satisfies the condition $0\leq v(x)\leq 1$ for almost all $x\in \mathcal{D}$. \textbf{Section 3} focuses on the finite element setup and discretization of the minimization problem. \textbf{Section 4} is divided into several parts, where the minimization algorithm is developed, residual-based estimates are derived, and the adaptive algorithm, which utilizes these estimates to improve the mesh, is presented. Additionally, convergence results are proven in Theorems \ref{Cngwithtol}--\ref{Cngwithouttol}, which follow the adaptive algorithms. Finally, \textbf{Section 5} concludes with numerical experiments that verify the analysis, followed by a summary of the findings.
\section{Elementary assets}
This section contains a few important aspects for our future investigation.
\subsection{Function spaces, norms and basic definitions}
In this work, we employ the usual notation from Lebesgue and Sobolev space theory \cite{adams1975}. The Sobolev space of order $(m,p)$, abbreviated by $\mbox{W}^m_p(\mathcal{D})$, is defined by  
\begin{equation*}
\mbox{W}^m_p(\mathcal{D})~:=\big\{\varphi \in \mbox{L}^p(\mathcal{D}):\quad D^l \varphi \in \mbox{L}^p(\mathcal{D}),\, |l|\leq m\big\}.
\end{equation*}
for $1\leq p \leq \infty$, equipped with inner product and the norm
\begin{equation*}
(\varphi, \psi)_{m,p,\mathcal{D}}~=~  \sum_{l \leq m}\int_\mathcal{D}  D^l \varphi \cdot D^l \psi\ , dx,   \quad 
\text{and} \quad  
\|\varphi\|_{\mbox{W}^m_p(\mathcal{D})}~=~ \Big(\sum_{|l| \leq m}\int_\mathcal{D}  |D^l \varphi|^p\, dx\Big)^{1/p},   
\end{equation*}
respectively. For $p=\infty$, the norm is given by
\begin{equation*}
\|\varphi\|_{\mbox{W}^m_\infty(\mathcal{D})}~=~ \esssup_{|l| \leq m}\|D^l \varphi\|_{L^\infty(\mathcal{D})}.   
\end{equation*}
Moreover, we write $\mbox{W}^m_2(\mathcal{D}):=\mbox{H}^m(\mathcal{D})$, for $p=2$, and the norm is denoted by $\|\cdot\|_{\mbox{H}^m(\mathcal{D})}$. Further, for $p=2$ and $m=1$, we now set $\mathbb{V}=\mbox{H}^1(\mathcal{D})$. We now introduce the function spaces for fixed $t=t_j$, which are crucial in our further study. Define
\begin{subequations}
\begin{eqnarray}
&& \mathbb{V}_d:=\big\{\varphi \in \mathbb{V}\,|~~\varphi=0~~ \text{on}~~ \Omega_D \big\},   \\
&& \mathbb{V}_c:=\big\{\varphi \in \mathbb{V}\,|~~ \varphi=0 ~~ \text{on}~~ CR(t_{j-1}) \big\},\\
\text{and} \nonumber\\
&& \mathbb{V}_f:=\big\{\varphi \in \mathbb{V}\,|~~ \varphi=f(t_j)~~ \text{on}~~ \Omega_D \big\},
\end{eqnarray}
\end{subequations}
respectively. 

Moreover, let $\mathscr{H}^d$ and $\mathscr{L}^d$ be the $d$-dimensional Hausdorff and Lebesgue measures, respectively. 
\begin{definition}[Total variation of $w$] Given a function $w \in L^1(\mathcal{D})$, then the total variation of $w$ in $\mathcal{D}$ is defined as 
\begin{align*}
\mathbb{V}(w, \mathcal{D}):=~ \sup\Big\{\int_{\mathcal{D}}w\,\div(\varphi)\,dx~:~~ \varphi \in C^1_0(\mathcal{D}; \mathbb{R}^d),\, |\varphi|\leq 1\Big\}.
\end{align*}
\end{definition}
\begin{definition}[Bounded variation]
The space of functions with bounded variation ($\mbox{BV}$ functions) can be defined as 
\begin{align*}
\mbox{BV}(\mathcal{D}):= \big\{w\in \mbox{L}^1(\mathcal{D})~:~~ \mathbb{V}(w, \mathcal{D})<+\infty \big\}. 
\end{align*}
\end{definition}
A function with bounded variation may have discontinuities that are represented in its distributional gradient $Dw$. This may be disassembled as
\begin{align*}
Dw~=~ \nabla w \mathscr{L}^d+ (w^+(x)-w^-(x)) \otimes \nu_{w}(x) \mathscr{H}^{d-1}\,\oldfactorial J(w)+D^cw,
\end{align*}
where $\nabla w$ and $J(w)$ denote the approximate gradient of $w$ and the jump set of $w$, while the $w^\pm$  denotes the inner and outer traces of $w$ on $J(w)$ with respect to $\nu_w$ (the unit normal vector to $J(w)$). Further, $D^c(w)$ represents the Cantor part of the derivative. For a more understanding of these notations, we refer to \cite{Ambrosio2002}. 
\begin{definition}[Special bounded variation] Let $SBV(\mathcal{D})$ be the collection of all bounded variation functions $w$ such that the Cantor part of the derivative of $w$ is zero, and defined as 
\begin{align*}
\mbox{SBV}(\mathcal{D}):=\big\{w\in \mbox{BV}(\mathcal{D})~:~~ D^c(w)=0\big\}.
\end{align*}   
\end{definition}
At this point, we are in the position to define the model problem.
\subsection{A mathematical framework} 

For $\kappa >0$ and $\epsilon >0$, we define the regularized elastic energy  $E \colon \mbox{H}^1(\mathcal{D}; \mathbb{R})\times \mbox{H}^1(\mathcal{D}; [0, 1]) \to \mathbb{R}^+_\infty$, where $\mathbb{R}^{+}_{\infty}=\mathbb{R}\cup \{+\infty\}$, and the regularized crack-surface energy $\mathcal{H} \colon \mbox{H}^1(\mathcal{D}; [0, 1]) \to \mathbb{R}$, respectively, as
\begin{subequations}\label{minmization}
\begin{align}
&E(u, \, v):= \frac{1}{2} \int_{\mathcal{D}} \dfrac{|T|^2}{\left(  1 + \beta^{\alpha}\,|T|^{2\alpha} \right)^{1/\alpha}}\; dx, \\ \text{and}~~~& \nonumber\\
 &\mathcal{H}(v):=\dfrac{1}{c_w}\int_{\mathcal{D}} \left[ \frac{(1-v)}{\epsilon} + \epsilon \; |\nabla v|^2  \right]  \; dx,
\end{align}
\end{subequations}
respectively, where $T:= T(u,v)= \sqrt{ (( 1- \kappa) v^2 + \kappa  )} \;  \nabla u$.
Then, the total energy is given by the functional  $\mathcal{J}_{\epsilon}: \mbox{H}^1(\mathcal{D}; \mathbb{R})\times \mbox{H}^1(\mathcal{D}; [0, 1]) \longmapsto \mathbb{R}$ such that
\begin{align}
\mathcal{J}_{\epsilon}(u, \, v) &:=   E(u, \, v) + \lambda_c \;  \mathcal{H}(v) \notag \\
&= \frac{1}{2} \int_{\mathcal{D}} \dfrac{|T|^2}{\left(  1 + \beta^{\alpha}\, | T|^{2\alpha} \right)^{1/\alpha}}\; dx +\dfrac{\lambda_c}{c_w}\int_{\mathcal{D}} \left[ \frac{(1-v)}{\epsilon} + \epsilon \; |\nabla v|^2  \right]  \; dx. \label{reg:energy}
\end{align}
Here $v \in \mbox{H}^{1} \left(\mathcal{D}; [0, \; 1] \right) $ represents a scalar \textit{phase-field function}, while the regularization parameter for the bulk energy term is considered to be 
$\kappa \ll 1$. Moreover, the normalization constant $c_w$  is chosen to be  $8/3$ for the \textit{damage field energy}.
Note that the phase-field variable $v \approx 0$ indicates an approximation to the crack-set, while $v=1$ represents the non-fractured zone, if
$v > 0$. However, the phase-field variable $v \in (0,1)$ depicted the critical length of the diffusive zone. For $\kappa\approx 0$, the bulk energy disappears, then the the crack energy only takes into account while $v=0$. In contrast, if $v = 1$, the fracture energy is zero, hence only the bulk energy needs to be taken into account. Therefore,  we interpolate nonzero fracture and bulk energies in the diffusive zone.
For the notational simplicity for boundary conditions, we define the following set 
\begin{equation} 
\mathrm{Q}(f(t)) :=\left\{ u \in \mbox{SBV}(\Omega) ~\colon~~ u|_{\Omega_D} = f(t) \right\},
\end{equation}
while the varying load $f\in \mbox{L}^\infty([0,T];\mbox{W}^1_{\infty}(\mathcal{D}))\cap \mbox{W}^1_1([0,T];\mbox{H}^1(\mathcal{D}))$ is applied on a open subset $\Omega_D\subset \mathcal{D}$. Consider the set of time points $\mathcal{S}_{t}=\{t_0,\, t_1,\, t_2,\,\ldots,\, t_{N_T}\}$ such that $0=t_0< t_1< t_2<t_3<\ldots< t_{N_T}=T$ with $k=\max\{k_j\,|\; k_j=t_j-t_{j-1},\; j=1,\,2,\ldots,\, N_T\}$ be the time points in a time domain $[0, T]$. So finally the quasi-static minimization problem reads as: 
 At time $t=t_0$, given $v(t_0, \, x) =1, \; \forall x \in \Omega$, find 
\begin{equation}
u_{\epsilon}(t_j, \, x) = \argmin \left\{E_{\epsilon}(\hat{u}, v(x)=1) \colon \hat{u} \in \mathbb{V}, \; \hat{ u}(x)=f(t_0) \;  \forall x \in \Omega_D \right\}.
\end{equation}
At the subsequent times $t=t_j$, $k=1, \ldots, m$, find $\left( u_{\epsilon}(t_j, \, x), \, v_{\epsilon}(t_j, \, x) \right)$ satisfying 
\begin{align}
\left( u_{\epsilon}(t_j), \, v_\epsilon(t_j, \, x) \right) = \argmin \big\{ &E_{\epsilon}(\hat{u},  \, \hat{v}) \colon \hat{u} \in \mathbb{V}, \; \hat{ u}(x)=f(t_j) \;  \forall x \in \Omega_D; \notag \\
&\hat{v} \in  \mathbb{V}, \, \hat{v} \leq v(t_{j-1}) \; \forall x \in \Omega \big\}. \label{eq_min_con}
\end{align}
The last term in the above equation \eqref{eq_min_con} imposes the crack irreversibility  \cite{Giacomini2005} and the condition is where we only allow the crack to propagate but not bonding or without self healing. At a fixed time $t=t_j$, the pair $(u_{\epsilon}(t_j), v_{\epsilon}(t_j))$ is an approximation for the Airy stress function $u(t_j)$ and the crack set $\Gamma(t_j)$ with
$u_{\epsilon}(t_j) \to u(t_j) \; \mbox{in} \; \mbox{L}^{1}(\Omega), \; \mbox{as} \; v \to 0.$

Since, the maximum principle concept is applicable to the phase-field variable. It is observed that the non-increasing character of the phase-field variable is a direct result of the truncation argument, making the following assertion easy to prove.
\begin{proposition}
The phase-field variable satisfy the maximum principle as 
\begin{equation}
0 \leq  v_{\epsilon}( x, \, t_j) \leq v_{\epsilon}(x, \, t_{j-1}) \quad \forall \; x \in \Omega, \quad \mbox{and} \quad \forall \; j \geq 1.
\end{equation}  
\end{proposition}
The quasi-static formulation can formulated as a successive global minimization  \cite{FrancfortBourdi2000} for the equilibrium solution as:
\begin{subequations}\label{eq:minmization}
\begin{align}
u(t_j) := \underset{v \in \mathcal{D}(f(t_j))}{\argmin} E_B(v) + E_S(J(v) + \Gamma(t_{k-1})), \quad \text{and} \quad 
\Gamma(t_j) := J(u(t_j)) \cup \Gamma(t_{k-1}).
\end{align}
\end{subequations}
Note that $(u(t_j), \Gamma(t_j))$ must satisfy the following crack irreversibility and global stability conditions.

\noindent 
\textbf{(i) \tt Crack irreversibility criteria:} The quasi-static evolution $t \mapsto \Gamma(t)$ requires that: 
\begin{equation}
\Gamma(t) \; \mbox{is increasing in time, i.e.} \; \Gamma(t_{j-1}) \subseteq \Gamma(t_j), \quad  \forall \; 0 \leq t_{j-1} \leq t_j \leq T
\end{equation}
The above condition can also be written as 
\begin{equation}\label{crack-irr-1}
 \partial_t v(x, t) \leq 0
 \end{equation}
  enforces the \textit{crack-irreversibility}. But, one can also impose the equality constraint as: If at time $t = t_{j-1}, \; j \in \left\{ 1, 2, \ldots, N_T \right\}$, define a set, for any $0 < c_{j} < 1$, 
\begin{equation}
v(t_{j-1}):= \left\{ x \in \overline{\mathcal{D}} \; \colon \; v_{\epsilon} (x,\, t_{j-1}) < c_{j} \right \},
\end{equation}
We note that the set $v(t_{j})$ is nonempty. Then the crack-irreversibility condition \eqref{crack-irr-1} be written as 
\begin{equation}
v_{\epsilon}(x, \, t_j) =0 \quad \forall \; x \in v(t_{j-1}) \; \mbox{and} \; \forall \; j \; \mbox{such that} \; k \leq j \leq N.
\end{equation} 
This simplifies minimization, as setting $v_{\epsilon}(x, \, t_j) $ to $0$ implies that the process be linear rather than nonlinear.  This clearly indicates that the point $x$ is on the crack path and continues to be on the crack-path until the end of the simulation. But the difficulty is that we may see \textit{crack widening} effect. This is easy to see, as setting $v_{\epsilon}(x, \, t_{j-1})$ to zero causes $v_{\epsilon}$ to remain below the value $c_{j}$ at points in a neighborhood of $x$, which will then be set to zero at $t=t_j$. One way to avoid the \textit{crack widening} effect is to choose $c_{j}$ less than the smallest mesh size. \smallskip

\noindent
\textbf{(ii) \tt Global stability criteria:} The sequence $(u(t_j), \Gamma(t_j)),\, j=1,\, \ldots,\, N_T$ for all $v\in \mathcal{Q}(f(t_{j-1}))$ and $\Gamma(t_j) \subset \Gamma$, satisfies 
\begin{align*}
E(u(t_j), \Gamma(t_j)) \leq E(v, \Gamma).
\end{align*}
For the problem formulated within non-linear elasticity, the existence of a solution to the Francfort-Marigo model has been shown in \cite{DeGiorgi1989,FrancfortLarsen2003}. But to the authors' best knowledge there is no such proof of existence of solutions available for the discrete-time quasi-static evolution model proposed within the setting of nonlinear limiting strain models introduced in \cite{yoon2021quasi}. 

Now we rewrite the functional $\mathcal{J}$ from $ \mathbb{V}_g\times \mathbb{V}_c$ to $\mathbb{R}^{+}_{\infty}$ by considering the constants $\kappa$, $\epsilon$, and $\lambda_c$ as fixed parameters:
\begin{align}
\mathcal{J}(u, \, v)~=~ \int_{\mathcal{D}}\Big[\dfrac{|T|^2}{2\left(  1 + \beta^{\alpha}\, |T|^{2\alpha} \right)^{1/\alpha}} + \rho\, |\nabla v|^2 + \delta \,(1-v)  \Big] \; dx, \label{minmizationfunc1.7}
\end{align}
with $\rho=\dfrac{\lambda_c\, \epsilon}{c_w }$ and  $\delta =\dfrac{\lambda_c}{c_w \,\epsilon}$, respectively. 

The following proposition demonstrates that the functional $\mathcal{J}$ is Fréchet differentiable at any $(u, v)\in \mathbb{V}\times \mathbb{V}_{T}$, where $\mathbb{V}_{T}=\mbox{H}^1(\mathcal{D})\cap L^\infty(\mathcal{D})$. 
\begin{proposition}[Fr\'{e}chet-differentiablity of $\mathcal{J}$]
For any $(u, v)\in \mathbb{V}\times \mathbb{V}_{T}$, then the functional  $\mathcal{J}$ from $\mathbb{V}_f\times \mathbb{V}_c $ to $\mathbb{R}^{+}_{\infty}$  is Fr\'{e}chet-differentiable in $\mathbb{V}\times \mathbb{V}_{T}$. 
\end{proposition}
\begin{proof} 
To examine the differentiability of $\mathcal{J}$ at $(u, v)\in \mathbb{V}\times \mathbb{V}_{T}$, we first compute the directional derivative in the direction $(\psi, \varphi)\in \mathbb{V}\times \mathbb{V}_{T}$, which yields
\begin{align}
\mathcal{J}^\prime(u, v; \psi, \varphi)=~ \mathcal{A}(v; u, \psi)+ \mathcal{B}(u; v, \varphi), \label{3.18jprim}
\end{align}
where
\begin{subequations}
\begin{align}
&\mathcal{A}(v; u, \psi)=~ \int_{\mathcal{D}}\Big[\dfrac{\big((1-\kappa)v^2+\kappa \big)}{\left(1 + \beta^{\alpha}\, |T|^{2\alpha}\right)^{\frac{1}{\alpha}+1}}\; \nabla u \cdot \nabla \psi \Big] \, dx, \label{ddofJA}\\
\text{and}~~~~& \nonumber\\
&\mathcal{B}(u; v, \varphi)=~ \int_{\mathcal{D}}\Big[2\,\rho\, \nabla v \cdot \nabla \varphi - \delta \, \varphi + \dfrac{(1-\kappa)\,|\nabla u|^2}{\left( 1 + \beta^{\alpha}\, |T|^{2\alpha}\right)^{\frac{1}{\alpha}+1}}\;v \, \varphi \Big] \, dx,  \label{ddofJb}
\end{align}   
\end{subequations}
respectively. To establish Fréchet differentiability, we need to verify that the remainder term $\frac{|\mathcal{R}(u, v; \psi, \varphi)|}{~\|\psi\|_{\mathbb{V}}+\|\varphi\|_{\mathbb{V}}+\|\varphi\|_{L^\infty(\mathcal{D})}}$ tends to zero at a rate consistent with $\|\psi\|_{\mathbb{V}}+\|\varphi\|_{\mathbb{V}}+\|\varphi\|_{L^\infty(\mathcal{D})}$ approaching zero. To this end, we will compute the remainder term $\mathcal{R}$ using the following definition
\begin{align}
&\big|\mathcal{R}(u, v; \psi, \varphi)\big|=~ \big|\mathcal{J}(u+\psi, v+\varphi)-\mathcal{J}(u, v)- \mathcal{J}^\prime(u, v; \psi, \varphi)\big|\nonumber\\
&~~~~~~\leq~\int_{\mathcal{D}}\rho\, |\nabla \varphi|^2\, dx+\Big|\int_{\mathcal{D}}\Big[\dfrac{((1-\kappa)(v+\varphi)^2+\kappa)|\nabla (u+\psi)|^2}{2\,[1 + \beta^{\alpha} ((1-\kappa)(v+\varphi)^2+\kappa)^\alpha|\nabla (u+\psi)|^{2\alpha} ]^{1/\alpha}}\nonumber\\
&~~~~~~~~~- \dfrac{((1-\kappa)v^2+\kappa)|\nabla u|^2}{2\,[1 + \beta^{\alpha} ((1-\kappa)v^2+\kappa)^\alpha\,|\nabla u|^{2\alpha} ]^{1/\alpha}} 
-\dfrac{\big((1-\kappa)v^2+\kappa \big)\; \nabla u \cdot \nabla \psi\,}{[ 1 + \beta^{\alpha}\,\big((1-\kappa)v^2+\kappa )\big)^\alpha\,| \nabla u|^{2\alpha}]^{\frac{1}{\alpha}+1}} \nonumber\\
&~~~~~~~~~ -\dfrac{(1-\kappa)\,|\nabla u|^2\;v \, \varphi}{[ 1 + \beta^{\alpha}\,\big((1-\kappa)v^2+\kappa )\big)^\alpha\, | \nabla u|^{2\alpha}]^{\frac{1}{\alpha}+1}} \Big]\,dx\Big| \nonumber\\
&~~~~~~:=~I_1+I_2. \label{3.20resbdd}
\end{align}
We proceed to evaluate each term individually. First, we estimate the term $I_1$ as
\begin{align*}
I_1\leq~ \rho\,\|\varphi\|^2_{\mathbb{V}}.
\end{align*}
Moving on to the next term using the fact $[ 1 + \beta^{\alpha}\,\big((1-\kappa)v^2+\kappa )\big)^\alpha | \nabla u|^{2\alpha}]^{\frac{1}{\alpha}+1}\geq [ 1 + \beta^{\alpha}\,\big((1-\kappa)v^2+\kappa )\big)^\alpha | \nabla u|^{2\alpha}]^{\frac{1}{\alpha}}$ and $1/[1 + \beta^{\alpha} ((1-\kappa)(v+\varphi)^2+\kappa)^\alpha\, |\nabla (u+\psi)|^{2\alpha} ]^{1/\alpha}\leq 1/[ 1 + \beta^{\alpha}\,\big((1-\kappa)v^2+\kappa )\big)^\alpha | \nabla u|^{2\alpha}]^{\frac{1}{\alpha}+1}$, we can estimate $I_2$ as
\begin{align*}
I_2&\leq~ \Big[ |1-\kappa|\,\|\varphi\|^2_{L^2(\mathcal{D})}\,\|\nabla u\|^2_{L^2(\mathcal{D})}+\|(1-\kappa)\,v^2+\kappa\|_{L^\infty(\mathcal{D})}\,\|\nabla \psi\|^2_{L^2(\mathcal{D})}+|1-\kappa|\,\|\varphi\|^2_{L^\infty(\mathcal{D})}\\
&~~\times \|\nabla \psi\|^2_{L^2(\mathcal{D})}+2\,|1-\kappa|\,\|v\|_{L^\infty(\mathcal{D})}\,\|\varphi\|_{L^\infty(\mathcal{D})}\,\|\nabla \psi\|^2_{L^2(\mathcal{D})}+4\,|1-\kappa|\,\|v\|_{L^\infty(\mathcal{D})}\,\,\|\nabla u\|_{L^2(\mathcal{D})}\\
&~~\times\|\varphi\|_{L^\infty(\mathcal{D})}\,\|\nabla \psi\|_{L^2(\mathcal{D})}+2\,|1-\kappa|\, \|\nabla u\|_{L^2(\mathcal{D})} \|\varphi\|^2_{L^\infty(\mathcal{D})}\,\|\nabla \psi\|_{L^2(\mathcal{D})}
\Big]{\Big/} [1 + \beta^{\alpha} \\ 
&~~\times \big((1-\kappa)v^2+\kappa )\big)^\alpha \| \nabla u\|^{2\alpha}]^{\frac{1}{\alpha}+1}.
\end{align*}
By substituting the estimates for $I_1$ and $I_2$ into \eqref{3.20resbdd}, we obtain $\frac{|\mathcal{R}(u, v; \psi, \varphi)|}{~\|\psi\|_{\mathbb{V}}+\|\varphi\|_{\mathbb{V}}+\|\varphi\|_{L^\infty(\mathcal{D})}}\longrightarrow 0$ as $\|\psi\|_{\mathbb{V}}+\|\varphi\|_{\mathbb{V}}+\|\varphi\|_{L^\infty(\mathcal{D})}\longrightarrow 0$. This concludes the proof of the lemma.
\end{proof}
It is observed that, the functional $\mathcal{J}$ lacks G$\hat{a}$teaux-differentiability over the entire $\mathbb{V}\times \mathbb{V}$. To circumvent this, we restrict our analysis to critical points within the subspace $\mathbb{V}_c^\infty=\mathbb{V}_c\cap \mbox{L}^\infty(\mathcal{D})$ of $\mathbb{V}_c$.
\begin{definition}[Critical point of $\mathcal{J}$] \label{crtitptsdef}
A pair $(u, v) \in \mathbb{V}_f \times \mathbb{V}_c^\infty$ is said to be a critical point of  $\mathcal{J}$ if  $(u, v)$ satisfy the condition $\mathcal{J}^{\prime}(u,v; \psi, \phi)=0$, $\forall \, (\psi, \phi) \in  \mathbb{V}_d \times \mathbb{V}_c^{\infty}$.
\end{definition}
It is well-known that any local minimizer pair $(u, v)$ of the functional $\mathcal{J}$ must satisfy the required condition $0\leq v(x)\leq 1$, which may be simply proven by using the truncation argument. So we omit the details of the proof, however, the statement is precisely laid out in the following proposition.
\begin{proposition} \label{prop3.2:v}
Let $(u, v)\in \mathbb{V}_f \times \mathbb{V}_c^\infty$ is a critical point of  $\mathcal{J}$, then $v \in \mathbb{V}_c^\infty$ satisfy the condition $0 \leq v(x) \leq 1$ for {\it a.e.} $x\in \mathcal{D}$.
\end{proposition}
\section{Finite element analysis}
The subsequent section describes the finite element setup for our model problem.
\subsection{Spatial discretization}
Let $\mathscr{T}_{h}$  be a family of regular simplicial triangulations of the domain $\bar{\mathcal{D}}$ such that
\begin{enumerate}
\item the boundary $\Gamma$ exactly represented by the boundaries of the triangles;
 for any two distinct triangles $\tau_i,\, \tau_j \in \mathscr{T}_{h},\, i \neq j$, then their intersection is either empty, a vertex, an edge, or a $k-$dimensional face, where $0\leq k \leq d-1$.
\item For each element $\tau \in \mathscr{T}_h$, we define its diameter as $h_\tau=diam(\tau)$ and the mesh size  $h=max\{h_\tau: \tau \in \mathscr{T}_h\}$. Further, 
the collection of all edges/faces are denoted by 
 $\mathscr{E}_{h}:=\mathscr{E}_{int, h} \cup \mathscr{E}_{bd,h},$ where the collection of all interior and the boundary edges represented by $\mathscr{E}_{int, h}$ and $\mathscr{E}_{bd,h}$, and are given by 
 $\mathscr{E}_{int, h}:= \{e_j \in \mathscr{E}_{h} \backslash (\bar{\Omega}_{D,h}\cup \partial \Omega_{N,h})\}$,
 $\mathscr{E}_{N, h}:= \{e_j \in \mathscr{E}_{h}:~~ e_j \subset \partial \Omega_N \},$ and 
$\mathscr{E}_{D, h}:= \{e_j \in \mathscr{E}_{h}:~~e_j \subset \bar{\Omega}_D \},$ respectively.
\item \texttt{Mesh shape regularity property.} We assume that the triangulation satisfies the shape-regularity condition: $\sup_{\tau \in \mathscr{T}_h} \frac{h_\tau}{\varrho_\tau} \leq c_{\varrho}$, where $\varrho_\tau$
denotes the largest diameter of the inscribed $d$-dimensional ball in $\tau$, and $c_{\varrho}$ is a positive constant.
\item \texttt{Reference element and affine map.} For each $\tau \in \mathscr{T}_h$, there exists an affine mapping $F:\hat{\tau}\mapsto \tau$, where  $\hat{\tau}$ is the reference simplex and defined as 
$$\hat{\tau}:= \big\{\hat{\mathtt{x}}~:~~ \hat{x}_j >0, \; \forall\, j\in [1:d], \quad 0<\sum_{j=1}^d \hat{x}_j<1 \big\}.$$ 
\item \texttt{Index set and basis functions.} We define the index set $\mathcal{N}_h \subset \mathbb{N}$ for the vertices of $\mathscr{T}_h$. We denote the basis function by $\xi_i$, $i \in \mathcal{N}_h$, such that
 $\xi_i$'s are continuous piecewise linear functions, and 
  $\xi_i(x_j)=\delta_{ij}$,  where $x_j$ represents the position of a vertex and $j\in \mathcal{N}_h$. 
\item \texttt{Useful notation.} We introduce the following notation which is useful in our analysis.
\begin{itemize}
\item[i.] Set $\mathcal{N}_{\Omega, h}=\big\{j \in \mathcal{N}_h~| ~~~ x_j \in \bar{\Omega} \big\},$ and
 $\omega_j$ is the closure of the union of elements $\tau_j \in \mathscr{T}_h$ that have $x_j$ as the position of a vertex, for $j\in \mathcal{N}_h$. That is, $\omega_j:= supp(\xi_j)$.
\item[ii.] Next, we set $h_{e_j}:=diam(e_j),$
 $\omega_{\tau_j}:= \bigcup_{\overset{j\in \mathcal{N}_h}{x_j \in \bar{\tau}_j}} \omega_j,$ and
 $\omega_e:=\bigcup_{\overset{j\in \mathcal{N}_h}{x_j \in \bar{e}_j}}  \omega_j,$ respectively.
\end{itemize}
\end{enumerate}
Exploiting the above phenomenon,  we introduce the finite element spaces $\mathbb{V}_h$, $\mathbb{V}_{d,h}$, and $\mathbb{V}^n_{f, h}$, providing a finite-dimensional framework, as 
\begin{align*}
& \mathbb{V}_h:=\big\{\sum_{j\in \mathcal{N}_h} \lambda_j\xi_j: ~~ \lambda_j\in \mathbb{R} \big\},  \nonumber\\
&\mathbb{V}_{d,h}:=\big\{\sum_{j\in \mathcal{N}_h} \lambda_j\xi_j: ~~ \lambda_j\in \mathbb{R},~~ \lambda_j=0~~ \text{for all}~j \in \mathcal{N}_{d,h}\big\},~ \text{and at time } t=t_n\\
&\mathbb{V}^n_{f, h}:=\big\{\sum_{j\in \mathcal{N}_h} \lambda_j\xi_j: ~~ \lambda_j\in \mathbb{R},~~ \lambda_j=f(t_n, x_j)~~ \text{for all}~j \in \mathcal{N}_{d,h} \big\}, 
\end{align*}
respectively. For a given tolerance ${\it \Xi_{CR}}$, and for $\varphi_h \in \mathbb{V}_h$, define a discrete version of $CR(t_{n-1})$, by
\begin{eqnarray*}
\mathcal{E}_h^{CR}(t_{n-1}):=\big\{ e_j \in \mathscr{E}_h\,\colon ~~ \varphi_h(x, t_{n-1})\leq {\it \Xi_{CR}},~~ \text{for all}~~ x\in \bar{e}_j \big\},
\end{eqnarray*}
such that  $CR_h(t_{n-1}):=\bigcup_{e_j \in \mathscr{E}^{CR}_h} \bar{e}_{j}$. Hence, the finite element space $\mathbb{V}^n_{c,h}$ is defined, as
\begin{eqnarray*}
\mathbb{V}^n_{c,h}:= \big\{\varphi_h\in \mathbb{V}_h~:~~ \varphi_h(x)=0,~~ \text{for all}~~ x\in CR_h(t_{n-1})\big\}.
\end{eqnarray*}
For simplicity, we write $\mathbb{V}_{c,h}^n$ and $\mathbb{V}_{f,h}^n$ by $\mathbb{V}_{c,h}$ and $\mathbb{V}_{f,h}$, respectively.

We are now able to develop a discrete version of the minimizing function $\mathcal{J}$. To achieve this, we first introduce the standard nodal interpolation operator $\pi_h: C(\bar{\Omega}) \mapsto \mathbb{V}_h$ \cite[Sec. 3.3]{Brenner1994}. Subsequently, we employ a mass lumping approximation (as described in \cite[Ch. 11]{vthomee1984}) to obtain the discrete formulation:
\begin{align}
\mathcal{J}_h(u_h, \, v_h)~=~ \int_{\mathcal{D}}\Big[\dfrac{|T^\pi_h|^2}{2(  1 + \beta^{\alpha}\, | T^\pi_h|^{2\alpha})^{1/\alpha}} + \rho\, |\nabla v_h|^2 + \delta\;\pi_h(1-v_h)\Big] \; dx, \label{disminmizationfunc4.1}
\end{align} 
where $T^\pi_h:=T^\pi_h(u_h,v_h)= \sqrt{( ( 1- \kappa)\; \pi_h(v_h^2) + \kappa  )} \;  \nabla u_h$. In addition, we will use $T_h= \sqrt{( ( 1- \kappa)\;v_h^2 + \kappa  )} \;  \nabla u_h.$
Next, we compute the critical point of $\mathcal{J}_h$. For this, we first introduce
\begin{subequations}
\begin{align}   
& \mathcal{J}_h^{\prime}(u_h,v_h; \psi_h, \varphi_h)~:=~  \mathcal{A}_h(v_h; u_h, \psi_h)+   \mathcal{B}_h(u_h; v_h, \varphi_h), \\
\text{where}~~~& \nonumber\\
&\mathcal{A}_h(v_h; u_h, \psi_h)=~ \int_{\mathcal{D}}\Big[\dfrac{\big((1-\kappa)\pi_h(v^2_h)+\kappa \big)}{\left(  1 + \beta^{\alpha}\, |T^\pi_h|^{2\alpha}\right)^{\frac{1}{\alpha}+1}}\; \nabla u_h \cdot \nabla \psi_h \Big] \, dx, \label{ddofJAfem}\\
\text{and}~~~~& \nonumber\\
&\mathcal{B}_h(u_h; v_h, \varphi_h)=~ \int_{\mathcal{D}}\Big[2\,\rho\, \nabla v_h \cdot \nabla \varphi_h - \delta \,\pi_h(\varphi_h) \nonumber\\
&\hspace{3.5cm}+ \dfrac{(1-\kappa)\,|\nabla u_h|^2}{\left( 1 + \beta^{\alpha}\, | T^\pi_h|^{2\alpha}\right)^{\frac{1}{\alpha}+1}}\;\pi_h(v_h \, \varphi_h) \Big] \, dx,  \label{ddofJbfem}
\end{align}   
\end{subequations}
respectively.
\begin{definition}[Critical point of $\mathcal{J}_h$]\label{def4.1} A pair $(u_h, v_h)\in \mathbb{V}_{f,h} \times \mathbb{V}_{c,h}$ is said to be a critical point of  $\mathcal{J}_h$ if  $(u_h, v_h)$ satisfy the condition $\mathcal{J}_h^{\prime}(u_h,v_h; \psi_h, \phi_h)=0$, $\forall (\psi_h, \phi_h) \in  \mathbb{V}_{d,h} \times \mathbb{V}_{c,h}$.
\end{definition}
Assuming the hypothesis (cf., \cite{ciarletreviart1973,  Korotov2001, Strangfix}) for the stiffness matrix $\mathbb{M}_{|\mathcal{N}_h|\times |\mathcal{N}_h|}$ (say), let $a_{ij}\in \mathbb{M}$, $1\leq i,j\leq \mathcal{N}_h$, satisfies the following condition:
\begin{eqnarray}
a_{ij}~:=~\int_{\mathcal{D}} \nabla \xi_i \cdot \nabla \xi_j\; dx~\leq~ 0, ~~\text{for all}~~ i\neq j \in  \mathcal{N}_h. \label{aij4.2}
\end{eqnarray}
Taking into account of inequality  \eqref{aij4.2}, one can easily show that the phase field minimizers $ v_h\in \mathbb{V}_{c,h}$, which are also a critical point, that lies between $0$ and $1$ for all $x\in \mathcal{D}$. We describe it as follows.
\begin{proposition} \label{prop4.2:vh}
Let $(u_h, v_h)$ be a pair in $ \mathbb{V}_{f,h} \times \mathbb{V}_{c,h}$ such that $\mathcal{B}_h(u_h; v_h, \phi_h)=0,~ \forall v_h \in \mathbb{V}_{c,h}$, then $v_h$ satisfy the condition $0 \leq v_h(x) \leq 1$, for all $x\in \mathcal{D}$.
\end{proposition}
The approximate model problem can be written as follows: Find a pair $(u_h, v_h)\in \mathbb{V}_{f,h} \times \mathbb{V}_{c,h}$ at $t=t_n$ so that
\begin{equation}
\mathcal{J}_h(u_h, v_h)=\underset{(\bar{u}_h, \bar{v}_h)\in \mathbb{V}_{f,h} \times \mathbb{V}_{c,h}}{argmin}\mathcal{J}_h(\bar{u}_h, \bar{v}_h). \label{minag4.3}
\end{equation}
To be more specific,
\begin{equation}
(u_h, v_h)\in argmin\big\{\mathcal{J}_h(\bar{u}_h, \bar{v}_h)~:~~(\bar{u}_h, \bar{v}_h)\in \mathbb{V}_{f,h} \times \mathbb{V}_{c,h}\big\}. \label{minagjh4.4}
\end{equation}
The following section offers alternative adaptive approaches to optimize the model \eqref{minag4.3}.
\section{Adaptive mechanisms and optimization} \label{minalgo}
To minimize the functional $\mathcal{J}$ on the infinite-dimensional space $\mathbb{V}_f\times \mathbb{V}_c$ at $t=t_n$ within tolerances $\Xi_v$ and  $\Xi_{v_n}$, we utilize the following minimization algorithm as stated below. 
\begin{flushleft}
\Hrule{4pt}{2pt} 
\vspace{0.1cm}
\textbf{Minimization Algorithm} \vspace{0.1cm}
\Hrule{4pt}{2pt} \vspace{0.2cm}

\textbf{Step 1.} {\tt Initialize:}  $v_0:=v(t_{j-1})$ if $j>0$~ and~ $v_0=1$ if $j=0$  \vspace{0.2cm}

\textbf{Step 2.} {\tt Start time loop:} For $t_j$, $j=1,\,2,\ldots,\, N_T.$ \\
\hspace{2cm} Compute $u_n$ and $v_n$ at $t_j$, for $n=1,\,2,\ldots$ \\
\hspace{2cm} Compute $u_n=\underset{\bar{u}\in \mathbb{V}_{f}}{{argmin}} \{\mathcal{J}(\bar{u}, v_{n-1})\}$ \\

\hspace{2cm} Compute $v_{n}=\underset{\bar{v}\in \mathbb{V}_{c}}{{argmin}} \{\mathcal{J}(u_n, \bar{v})\}$
\\ \vspace{0.2cm}

\hspace{0.7cm} \textbf{Check:} If $\|v_{n}-v_{n-1}\|_{L^\infty(\mathcal{D})}\geq \Xi_{v_n}$\\

\hspace{3cm} Repeat \textbf{\tt Step 2} \\

\hspace{2.5cm} Else  \\

\hspace{3.0cm} \textbf{\tt Break}; \\

\hspace{2.5cm} End Else \\
\hspace{3cm} If $v_n \leq \Xi_v $\\
\hspace{3.5cm} $v_n=0$\\
\hspace{3.5cm} Else if $v_n\geq 1.0$\\
\hspace{4.0cm} $ v_n=1.0$\\
\hspace{3.5cm} End Else if\\
\hspace{3cm} End If \\
\hspace{2.0cm} End If \\ \vspace{0.2cm}

\textbf{Step 3.} Set $u_n=u(t_j)$ and $v_n=v(t_j)$. \\ \vspace{0.1cm}
\hspace{2 cm} Repeat the steps. \vspace{0.2cm}
\Hrule{4pt}{2pt}
\end{flushleft} \vspace{0.1cm}
\noindent
The convergence of our numerical solution is analysed in the following sections, with Theorem \ref{Cngwithtol} and Theorem\ref{Cngwithouttol} providing a summary of the key findings. The primary distinction between the following algorithms lies in the timing and location of mesh refinement. We propose two adaptive algorithms:
(1) \texttt{Post-optimization refinement technique} (Adaptive Algorithm-1) and;  
(2) \texttt{Iterative optimization refinement technique} (Adaptive Algorithm-2).
It is important to note that the \texttt{Adaptive Algorithm-1} refines the mesh post-minimization after the completion of the alternate Minimization Algorithm (Sec. \ref{minalgo}), whereas the \texttt{Adaptive Algorithm-2} incorporates mesh refinement at each stage of the alternate minimization procedure.

To make it easier to construct these algorithms, we use approximation features to provide residual-type indicators.
Let $\pi^d_h$ and $\pi^c_h$ be the {\it quasi interpolants} such that $\pi^d_h \phi \in \mathbb{V}_{d,h}$ and $\pi^c_h \phi \in \mathbb{V}_{c,h}$ for $\phi \in \mathbb{V}$ and defined as 
\begin{subequations} 
\begin{align}
&\big(\pi^d_h \phi \big)(x) := \sum_{j\in \mathcal{N}_h \setminus \mathcal{N}_{d,h}} \Big(\frac{1}{|M_j|}\int_{M_j} \phi \,dx \Big)(x), \label{app4.6aa}\\
\text{and}~~~~& \nonumber\\ 
&\big(\pi^c_h \phi \big)(x) := \sum_{\overset{j\in \mathcal{N}_h}{x_j \notin CR_h}} \Big(\frac{1}{|M_j|}\int_{M_j} \phi \,dx \Big)(x), \label{app4.6bb}
\end{align}
\end{subequations}
respectively, where $M_j$ denotes the maximal set contained in $\omega_j$ (cf., \cite{Verfurth1999}).
\begin{lemma}[\tt{Quasi-interpolation approximation properties}] \label{4.3lemmaapr}
Let $\pi^d_h$ and $\pi^c_h$ be the {\it quasi interpolants} defined as \eqref{app4.6aa}-\eqref{app4.6bb}. Then, there exist positive constant  $c_{41}$, $c_{42}$, $c_{43}$ and $c_{44}$, the constants may depend on the shape-regularity parameter of the mesh $\tau$ but not on the mesh size, such that, for all elements $\tau_j \in \mathscr{T}_h$  and edges $e_j \in \mathscr{E}_h$, $\forall\;  j\in \mathcal{N}_h$, for $\phi \in \mathbb{V}_d$ and  $m \in \{0, 1\}$,
\begin{subequations}
\begin{align}
& \|\phi-\pi^d_h \phi\|_{H^m(\tau_j)} \leq c_{41}\,h^{1-m}_{\tau_j} \|\nabla \phi\|_{L^2(\omega_{\tau_j})}, ~~ \|\phi-\pi^d_h \phi\|_{L^2(e_j)} \leq c_{42}\,h^{1/2}_{e_j} \|\nabla \phi\|_{L^2(\omega_{e_j})}, \label{4.10AAapr} \\
\text{and}&~\text{for } \phi \in \mathbb{V}_c,  \nonumber\\
& \|\phi-\pi^c_h \phi\|_{H^m(\tau_j)} \leq c_{43}\,h^{1-m}_{\tau_j} \|\nabla \phi\|_{L^2(\omega_{\tau_j})} ,  \|\phi-\pi^c_h \phi\|_{L^2(e_j)} \leq c_{44}\,h^{1/2}_{e_j} \|\nabla \phi\|_{L^2(\omega_{e_j})}. \label{4.10bbapr}
\end{align}
\end{subequations}
\end{lemma}
\begin{lemma}[\tt{Standard nodal interpolation approximation properties}] \label{4.3lemmaaprLinfty}
Let $\pi_h$ be the standard nodal interpolation operator as defined in  \cite[Sec 3.3]{Brenner1994}. For all $\phi \in W^m_\infty(\tau_j)$, $\forall\, \tau_j \in \mathscr{T}_h,$ $\forall\;  j\in \mathcal{N}_h$ and $m \in \{1, 2\}$, then there exist positive constant  $c_{45}$, such that
\begin{align}
& \|\phi-\pi_h \phi\|_{L^\infty(\tau_j)} \leq c_{45}\,h^m_{\tau_j} |\phi|_{W^m_\infty(\tau_j)}. \label{4.11appstd}
\end{align}
\end{lemma}
For the estimation, we define the notion of a jump for any $\phi_h \in \mathbb{V}_h$ across an internal element edge/face $e_{ij}\in \mathscr{E}_{int,h}$ shared by elements $\tau_i$ and $\tau_j$ ($i>j$), as follows:
\begin{align*}
\sjump{\nabla \phi_h}:=~|\nabla \phi_h|_{e_{ij}\cap \partial\tau_i}-|\nabla \phi_h|_{e_{ij}\cap \partial\tau_j},\quad \text{and} \quad \sjump{\nabla \phi_h}|_{e_{j}}=:~|\nabla \phi_h \cdot \textbf{n}|_{\partial \tau_{j}\cap \Gamma},\;\; e_j \in \partial \tau_j \cap \Gamma,
\end{align*}
where $\textbf{n}$ is the outer unit normal vector.

The following lemma constitutes the foundation for our convergence analysis in the following section.
\begin{lemma} \label{4.5jbddlemma}
For all $\psi_h \in \mathbb{V}_{d,h}$ and $\varphi_h \in \mathbb{V}_{c,h}$, we assume that $u_h\in \mathbb{V}_{f,h}$ and $v_h\in \mathbb{V}_{c,h}$ such that $\mathcal{J}^\prime(u_h, v_h; \psi_h, \varphi_h) =0 $. Then, there exists $c_{46}>0$ such that the following estimate holds,  for all $\psi \in \mathbb{V}_d$ and $\varphi \in \mathbb{V}_c$, 
\begin{align}
 |\mathcal{J}^\prime(u_h, v_h; \psi, \varphi) | \leq~ c_{46}\, \big\{\widetilde{\eta}_h\,\|\nabla \psi\|_{L^2(\mathcal{D})}+\widehat{\eta}_h\,\|\nabla \varphi\|_{L^2(\mathcal{D})}\big\},    \label{4.12lemmabdd}
\end{align}
where $\widetilde{\eta}_h$ and $\widehat{\eta}_h$ are defined as 
\begin{align}
\widetilde{\eta}_h:=~\Big(\sum_{\tau \in \mathscr{T}_h} |\widetilde{\eta}_{\tau}(u_h,v_h)|^2 \Big)^{\frac{1}{2}}, \quad \text{and} \quad \widehat{\eta}_h:=~\Big(\sum_{\tau \in \mathscr{T}_h} |\widehat{\eta}_{\tau}(u_h,v_h)|^2 \Big)^{\frac{1}{2}}
\end{align}
respectively, where 
\begin{subequations} 
\begin{align} 
|\widetilde{\eta}_{\tau}(u_h,v_h)|^2 &:=~\sup_{x\in \tau_i}|\nabla v_h|^4\times\int_{\tau_i} h_{\tau_i}^4\,\Big|\dfrac{(1-\kappa)\, \nabla u_h }{\big[ 1 + \beta^{\alpha}\,| T_h |^{2\alpha} \big]^{\frac{1}{\alpha}+1}} \Big|^2\,dx \nonumber\\ 
&~~~~+\int_{\tau_i} h_{\tau_i}^2\, \Big| \dfrac{2\,(\kappa-1)v_h \; \big(\nabla v_h\cdot \nabla u_h\big) [1-\alpha\, \beta^\alpha\,|T_h|^{2\alpha} ] }{\big[1 + \beta^{\alpha}\,\,|T_h|^{2\alpha} \big]^{\frac{1}{\alpha}+2}}  \Big|^2 dx \nonumber\\ 
&~~~~+\sum_{ e_i\in \partial \tau_i \cap (\mathscr{E}_{int, h} \cup \mathscr{E}_{N, h}) }\int_{e_i} h_{e_i}\Big|\dfrac{\big((1-\kappa)v_h^2+\kappa \big)\, \sjump{\nabla u_h}}{\big[ 1 + \beta^{\alpha}\,\,|T_h|^{2\alpha} \big]^{\frac{1}{\alpha}+1}}\Big|^2\,ds \label{esteta4.11aa}\\ 
\text{and}~~~~~~~~~~~& \nonumber \\ 
|\widehat{\eta}_{\tau}(u_h,v_h)|^2 &:=~ \sup_{x\in \tau_i}|\nabla v_h|^2\times\int_{\tau_i} h_{\tau_i}^4\,\Big|\dfrac{(1-\kappa)\,|\nabla u_h|^2}{\left[ 1 + \beta^{\alpha}\, | T_h|^{2\alpha}\right]^{\frac{1}{\alpha}+1}}\Big|^2\,dx  \nonumber\\ 
&~~~~+ \int_{\tau_i} h^2_{\tau_i}\,\Big|\dfrac{(1-\kappa)\,|\nabla u_h|^2\,v_h}{\left[ 1 + \beta^{\alpha}\,| T_h |^{2\alpha}\right]^{\frac{1}{\alpha}+1}}- \delta \Big|^2\,dx+ \sum_{ e_i\in \partial \tau_i\cap \mathscr{E}_h} \int_{e_i} \rho^2 \,h_{e_i}\,|\sjump{\nabla v_h}|^2\,ds. \label{esteta4.11bb}
\end{align}
\end{subequations}
Here $c_{46}=\max\{c_{47}, c_{48}\} $, where the constants $c_{47}$ and $ c_{48}$ are defined by $c_{47}=\max\{c_{41}, c_{42}, c_{45}\}$ and  $c_{48}=\max\{c_{43},\, 2\, c_{44}, c_{\varrho}\,c_{eq}\,c_{45}\}$, respectively. 
\end{lemma}
\begin{proof}
For every $\psi \in \mathbb{V}_d$ and $\varphi \in \mathbb{V}_c$, the equation \eqref{3.18jprim} implying that 
\begin{align}
\big| \mathcal{J}^\prime(u_h, v_h; \psi, \varphi)\big|\leq~ \big|\mathcal{A}(v_h; u_h, \psi)\big|+ \big| \mathcal{B}(u_h; v_h, \varphi)\big|, \label{483.18jprim}
\end{align}
To begin, we will establish individual bounds for the terms $\mathcal{A}$  and $\mathcal{B}$. We start by bounding the term $\mathcal{A}$ for any $\psi\in \mathbb{V}_d$ and $\psi_h \in \mathbb{V}_{d,h}$, yielding
\begin{align}
\big|\mathcal{A}(v_h; u_h, \psi)\big|=  \big|\mathcal{A}(v_h; u_h, \psi-\psi_h)\big|+\big|\mathcal{A}(v_h; u_h,\psi_h)-\mathcal{A}_h(v_h; u_h,\psi_h)\big|:=I_1+I_2. \label{4.16Abd}
\end{align}
Since $\mathcal{A}_h(v_h; u_h,\psi_h)=0$ for all $\psi_h \in \mathbb{V}_{d,h}$.  We initiate our analysis with the term $I_1$. Applying the Green's inequality, we obtain
\begin{align*}
I_1&= ~ \Big| \sum_{ \tau_i\in \mathscr{T}_h}\int_{\tau_i}\Big[\dfrac{\big((1-\kappa)v_h^2+\kappa \big)}{\big[ 1 + \beta^{\alpha}\, | T_h|^{2\alpha} \big]^{\frac{1}{\alpha}+1}}\; \nabla u_h \cdot \nabla (\psi-\psi_h) \Big] \, dx \Big|\\
&\leq \Big(\sum_{ \tau_i\in \mathscr{T}_h} \int_{\tau_i} \Big| \dfrac{-2(1-\kappa)v_h \; \big(\nabla v_h\cdot \nabla u_h\big) [1-\alpha\, \beta^\alpha\,|T_h|^{2\alpha} ] }{\big[1 + \beta^{\alpha}\,|T_h|^{2\alpha} \big]^{\frac{1}{\alpha}+2}}  \Big|^2 dx\Big)^{\frac{1}{2}}\times \|\psi-\psi_h\|_{L^2(\tau_i)} \\ 
&~~~~ +  \Big(\sum_{ \tau_i\in \mathscr{T}_h}\sum_{ e_i\in \partial \tau_i \cap (\mathscr{E}_{int, h} \cup \mathscr{E}_{N, h}) }\int_{e_i}\Big|\dfrac{\big((1-\kappa)v_h^2+\kappa \big)\, \sjump{\nabla u_h}}{\big[ 1 + \beta^{\alpha}\,|T_h|^{2\alpha} \big]^{\frac{1}{\alpha}+1}}\Big|^2\,ds\Big)^{\frac{1}{2}}\times \|\psi-\psi_h\|_{L^2(e_i)}.
\end{align*}
Setting $\psi_h=\pi_h^d \psi$, then using the approximation properties \ref{4.10AAapr} of Lemma \ref{4.3lemmaapr} with $m=0$ and the Cauchy-Schwarz inequality, we arrive at  
\begin{align}
I_1 &\leq c_{41}\, \Big(\sum_{ \tau_i\in \mathscr{T}_h} \int_{\tau_i} h_{\tau_i}^2\, \Big| \dfrac{-2(1-\kappa)v_h \; \big(\nabla v_h\cdot \nabla u_h\big) [1-\alpha\, \beta^\alpha\,|T_h|^{2\alpha} ] }{\big[1 + \beta^{\alpha}\,\,|T_h|^{2\alpha} \big]^{\frac{1}{\alpha}+2}}  \Big|^2 dx\Big)^{\frac{1}{2}} \times  \|\nabla \psi\|_{L^2(\omega_{\tau_i})} \nonumber\\
&~~ + c_{42} \,\Big(\sum_{ \tau_i\in \mathscr{T}_h}\sum_{ e_i\in \partial \tau_i \cap (\mathscr{E}_{int, h} \cup \mathscr{E}_{N, h}) }\int_{e_i}h_{e_i}\Big|\dfrac{\big((1-\kappa)v_h^2+\kappa \big)\, \sjump{\nabla u_h}}{\big[ 1 + \beta^{\alpha}\,|T_h|^{2\alpha} \big]^{\frac{1}{\alpha}+1}}\Big|^2\,ds\Big)^{\frac{1}{2}}\times \|\nabla \psi\|_{L^2(\omega_{e_i})}\nonumber\\
&\leq \max\{c_{41}, c_{42} \} \times \Bigg\{\sum_{ \tau_i\in \mathscr{T}_h} \int_{\tau_i} h_{\tau_i}^2\, \Big| \dfrac{-2(1-\kappa)v_h \; \big(\nabla v_h\cdot \nabla u_h\big) [1-\alpha\, \beta^\alpha\,|T_h|^{2\alpha} ] }{\big[1 + \beta^{\alpha}\,|T_h|^{2\alpha} \big]^{\frac{1}{\alpha}+2}}  \Big|^2 dx \nonumber\\
&~~~+ \sum_{ \tau_i\in \mathscr{T}_h}\sum_{ e_i\in \partial \tau_i \cap (\mathscr{E}_{int, h} \cup \mathscr{E}_{N, h}) }\int_{e_i}h_{e_i}\Big|\dfrac{\big((1-\kappa)v_h^2+\kappa \big)\, \sjump{\nabla u_h}}{\big[ 1 + \beta^{\alpha}\,|T_h|^{2\alpha} \big]^{\frac{1}{\alpha}+1}}\Big|^2\,ds \Bigg\}^{\frac{1}{2}} \nonumber \\
&~~~ \times \Bigg\{ \sum_{ \tau_i\in \mathscr{T}_h}\Big[\|\nabla \psi\|^2_{L^2(\omega_{\tau_i})}+\sum_{ e_i\in \mathscr{E}_{int, h} \cup \mathscr{E}_{N, h} }  \|\nabla \psi\|^2_{L^2(e_i)}\Big] \Bigg\}^{\frac{1}{2}} \nonumber\\
&\leq \max\{c_{41}, c_{42} \} \times \Bigg\{\sum_{\tau_i\in \mathscr{T}_h} \bigg[\int_{\tau_i} h_{\tau_i}^2\, \Big| \dfrac{2\,(\kappa-1)v_h \; \big(\nabla v_h\cdot \nabla u_h\big) [1-\alpha\, \beta^\alpha\,|T_h|^{2\alpha} ] }{\big[1 + \beta^{\alpha}\,|T_h|^{2\alpha} \big]^{\frac{1}{\alpha}+2}}  \Big|^2 dx \nonumber\\
&~~~+\sum_{ e_i\in \partial \tau_i \cap (\mathscr{E}_{int, h} \cup \mathscr{E}_{N, h}) }\int_{e_i} h_{e_i}\Big|\dfrac{\big((1-\kappa)v_h^2+\kappa \big)\, \sjump{\nabla u_h}}{\big[ 1 + \beta^{\alpha}\,|T_h|^{2\alpha} \big]^{\frac{1}{\alpha}+1}}\Big|^2\,ds \bigg] \Bigg\}^{\frac{1}{2}}\times  \|\nabla \psi\|_{L^2(\mathcal{D})}.\label{4.17I1bdd}
\end{align}
Next, we bound the term $I_2$ utilizing the operator property $0\leq \pi_h(v_h^2)\leq v_h^2$, which leads to 
\begin{align}
 I_2 &= \Big| \sum_{ \tau_i\in \mathscr{T}_h}\int_{\tau_i}\Big[\dfrac{\big((1-\kappa)v_h^2+\kappa \big)}{\big[ 1 + \beta^{\alpha} | T_h|^{2\alpha} \big]^{\frac{1}{\alpha}+1}} -\dfrac{\big((1-\kappa)\pi_h(v^2_h)+\kappa \big)}{\left(  1 + \beta^{\alpha}\, |T_h^\pi|^{2\alpha}\right)^{\frac{1}{\alpha}+1}}\Big]\, \nabla u_h \cdot \nabla \psi_h  \, dx \Big| \nonumber\\
 &\leq \Big| \sum_{ \tau_i\in \mathscr{T}_h}\int_{\tau_i} \dfrac{(1-\kappa)(v_h^2-\pi_h(v^2_h))}{\big[ 1 + \beta^{\alpha} \, | T_h |^{2\alpha} \big]^{\frac{1}{\alpha}+1}} \, \nabla u_h \cdot \nabla \psi_h  \, dx \Big| \nonumber\\
 & \leq \Big\{ \sum_{ \tau_i\in \mathscr{T}_h}\int_{\tau_i} \Big| \dfrac{(1-\kappa)\, \nabla u_h }{\big[ 1 + \beta^{\alpha}\, | T_h |^{2\alpha} \big]^{\frac{1}{\alpha}+1}} \Big|^2\,dx \Big\}^{\frac{1}{2}} \times \|v_h^2-\pi_h(v^2_h)\|_{L^\infty(\tau_i)}\; \|\nabla \psi_h \|_{L^2(\tau_i)}.\nonumber 
\end{align}
Applying the nodal approximation property [\ref{4.11appstd}, Lemma \ref{4.3lemmaaprLinfty}]  with $m=1$, to obtain 
\begin{align}
 I_2 & \leq c_{45}\,\Big\{ \sum_{ \tau_i\in \mathscr{T}_h}\sup_{x\in \tau_i}|\nabla v_h|^4\times\int_{\tau_i} h_{\tau_i}^4\Big|\dfrac{(1-\kappa)\, \nabla u_h }{\big[ 1 + \beta^{\alpha}\, |T_h |^{2\alpha} \big]^{\frac{1}{\alpha}+1}} \Big|^2\,dx \Big\}^{\frac{1}{2}}\, \|\nabla \psi_h \|_{L^2(\omega_{\tau_i})} \nonumber \\
 &\leq  c_{45}\,\Big\{ \sum_{ \tau_i\in \mathscr{T}_h}\sup_{x\in \tau_i}|\nabla v_h|^4\times\int_{\tau_i} h_{\tau_i}^4\,\Big|\dfrac{(1-\kappa)\, \nabla u_h }{\big[ 1 + \beta^{\alpha}\, |T_h |^{2\alpha} \big]^{\frac{1}{\alpha}+1}} \Big|^2\,dx \Big\}^{\frac{1}{2}}\, \|\nabla \psi \|_{L^2(\mathcal{D})}. \label{4.18I2}
\end{align}
Plugging the bounds for $I_1$ and $I_2$ from \eqref{4.17I1bdd} and \eqref{4.18I2} into \eqref{4.16Abd} yields
\begin{align}
\big|\mathcal{A}(v_h; u_h, \psi)\big|&= \max\{c_{41}, c_{42}, c_{45}\} \times \Big\{\sum_{\tau \in \mathscr{T}_h} |\widetilde{\eta}_{\tau}(u_h,v_h)|^2 \Big\}^{\frac{1}{2}}\times  \|\nabla \psi\|_{L^2(\mathcal{D})}.\label{4.19auhpsibdd}
\end{align}
We now derive an estimate for the term $|\mathcal{B}(u_h;v_h, \varphi)|$, for all $\varphi \in \mathbb{V}_c$. Note that for any $\varphi_h \in \mathbb{V}_{c,h}$, the term $|\mathcal{B}(u_h;v_h, \varphi)|$ can be bounded as follows 
\begin{align}
|\mathcal{B}(u_h;v_h, \varphi)|\leq~ |\mathcal{B}(u_h;v_h, \varphi-\varphi_h)|+|\mathcal{B}(u_h;v_h, \varphi_h)- \mathcal{B}_h(u_h;v_h, \varphi_h)|:=I_3+I_4. \label{4.20BB}
\end{align}
Since $\mathcal{B}_h(u_h;v_h, \varphi_h)=0,\, \forall \varphi_h \in \mathbb{V}_{c,h}$. 
We proceed to estimate the terms $I_3$ and $I_4$, separately. Beginning with $I_3$, and utilizing \eqref{ddofJb} by choosing an arbitrary $\varphi_h \in \mathbb{V}_{c,h}$, we derive the following bound
\begin{align}
I_3 &\leq \Big|\sum_{ \tau_i\in \mathscr{T}_h}\int_{\tau_i} \Big[\dfrac{(1-\kappa)\,|\nabla u_h|^2\;v_h}{\left[ 1 + \beta^{\alpha}\,|T_h|^{2\alpha}\right]^{\frac{1}{\alpha}+1}}- \delta \Big]\,(\varphi-\varphi_h) \, dx\Big|+\Big|\sum_{ \tau_i\in \mathscr{T}_h}\int_{\partial \tau_i} 2\,\rho\, \nabla v_h \cdot \textbf{n}\; (\varphi - \varphi_h) \, dx \Big| \nonumber\\ 
&\leq \Big\{\sum_{ \tau_i\in \mathscr{T}_h} \int_{\tau_i}\Big|\dfrac{(1-\kappa)\,|\nabla u_h|^2\,v_h}{\left[ 1 + \beta^{\alpha}\,|T_h|^{2\alpha}\right]^{\frac{1}{\alpha}+1}}-\delta \Big|^2\,dx\Big\}^{\frac{1}{2}}\|\varphi-\varphi_h\|_{L^2(\tau_i)}\nonumber\\
&~~~+2\,\rho\,\Big\{\sum_{ \tau_i\in \mathscr{T}_h} \sum_{ e_i\in \partial \tau_i\cap \mathscr{E}_h}\int_{e_i}  |\sjump{\nabla v_h}|^2\,ds\Big\}^{\frac{1}{2}}\|\varphi - \varphi_h\|_{L^2(e_i)}.\nonumber
\end{align}
Setting $\varphi_h=\pi_h^c \varphi$, then use of inequality [\ref{4.10bbapr}, Lemma \ref{4.3lemmaapr}] with $m=0$ and the Cauchy Schwarz inequality leads to 
\begin{align}
I_3&\leq \max\{c_{43},\, 2\, c_{44}\} \times \Big\{\sum_{ \tau_i\in \mathscr{T}_h} \Big[\int_{\tau_i} h^2_{\tau_i}\,\Big|\dfrac{(1-\kappa)\,|\nabla u_h|^2\,v_h}{\left[ 1 + \beta^{\alpha}\,|T_h |^{2\alpha}\right]^{\frac{1}{\alpha}+1}}-\delta \Big|^2\,dx\nonumber\\
&~~ + \sum_{ e_i\in \partial \tau_i\cap \mathscr{E}_h} \int_{e_i} \rho^2 \,h_{e_i}\,|\sjump{\nabla v_h}|^2\,ds \Big] \Big\}^{\frac{1}{2}}\times \Big\{\sum_{ \tau_i\in \mathscr{T}_h}\Big[ \|\nabla \varphi\|^2_{L^2(\omega_{\tau_i})}+ \sum_{ e_i\in \partial \tau_i \cap \mathscr{E}_h} \|\nabla \varphi\|^2_{L^2(\omega_{e_i})} \Big] \Big\}^{\frac{1}{2}}.\label{4.21BBbdd}
\end{align}
We now turn our attention to estimating the term $I_4$. To this end, we employ the standard nodal interpolation estimate as stated in Lemma \ref{4.3lemmaaprLinfty} with utilizing the operator property, and then by setting $\varphi_h=\pi_h^c \varphi$, we bound the term $I_4$ as 
\begin{align}
I_4&=~\Big|\sum_{ \tau_i\in \mathscr{T}_h}\int_{\tau_i} \Big[\Big(\dfrac{(1-\kappa)\,|\nabla u_h|^2}{\left[ 1 + \beta^{\alpha} |T_h|^{2\alpha}\right]^{\frac{1}{\alpha}+1}}\Big)\;(v_h\varphi_h -\pi_h(v_h \, \varphi_h))\Big]\, dx \,\Big| \nonumber\\
&\leq~ \Big\{\sum_{ \tau_i\in \mathscr{T}_h} \int_{\tau_i} h_{\tau_i}^d\,\Big|\dfrac{(1-\kappa)\,|\nabla u_h|^2}{\left[ 1 + \beta^{\alpha}\,| T_h|^{2\alpha}\right]^{\frac{1}{\alpha}+1}} \Big|^2\,dx \Big\}^{\frac{1}{2}}\, \|v_h\varphi_h -\pi_h(v_h \, \varphi_h)\|_{L^\infty(\tau_i)}  \nonumber \\
&\leq~  c_{45}\,\Big\{\sum_{ \tau_i\in \mathscr{T}_h} \int_{\tau_i} h_{\tau_i}^{d+4}\,\Big|\dfrac{(1-\kappa)\,|\nabla u_h|^2}{\left[ 1 + \beta^{\alpha}\,| T_h|^{2\alpha}\right]^{\frac{1}{\alpha}+1}} \Big|^2\,dx \Big\}^{\frac{1}{2}}\times\|v_h\varphi_h \|_{W^2_\infty(\tau_i)}.  \nonumber
\end{align}
By virtue of norm equivalence in finite-dimensional spaces and the shape regularity condition, it follows that
\begin{align}
I_4 &\leq~ c_{\varrho}\,c_{eq}\,c_{45}\,\Big\{\sum_{ \tau_i\in \mathscr{T}_h} 
\sup_{x\in \tau_i}|\nabla v_h|^2\times \int_{\tau_i} h_{\tau_i}^4\,\Big|\dfrac{(1-\kappa)\,|\nabla u_h|^2}{\left[ 1 + \beta^{\alpha}\, | T_h|^{2\alpha}\right]^{\frac{1}{\alpha}+1}}\Big|^2\,dx\Big\}^{\frac{1}{2}}\times \|\nabla \varphi\|_{L^2(\mathcal{D})},\label{4.22I4bdd}  
\end{align}
Here $c_\varrho$ and   $c_{eq}$ denote the shape regularity  and the norm equivalence constants, respectively. Combining the estimates for $I_3$ and $I_4$ from \eqref{4.21BBbdd} and \eqref{4.22I4bdd}, respectively, into \eqref{4.20BB} yields
\begin{align}
\big|\mathcal{B}(u_h;v_h,\varphi)\big|&\leq~ \max\{c_{43},\, 2\, c_{44}, c_{\varrho}\,c_{eq}\,c_{45}\} \times \Big\{ \sum_{ \tau_i\in \mathscr{T}_h} |\widehat{\eta}_{\tau}(u_h,v_h)|^2\Big\}^{\frac{1}{2}} \times  \|\nabla \varphi\|_{L^2(\mathcal{D})}. \label{4.23buhbdd}
\end{align}
Merging the derived estimates of $\big|\mathcal{A}(v_h;u_h,\psi)\big|$ and $\big|\mathcal{B}(u_h;v_h,\varphi)\big|$ from \eqref{4.19auhpsibdd} and \eqref{4.23buhbdd} with \eqref{3.18jprim}, and setting $c_{46}=\max\{c_{47}, c_{48}\} $, where $c_{47}=\max\{c_{41}, c_{42}, c_{45}\}$ and  $c_{48}=\max\{c_{43},\, 2\, c_{44}, c_{\varrho}\,c_{eq}\,c_{45}\}$, respectively, this leads to the desired inequality \eqref{4.12lemmabdd}. Therefore, the proof is completed.
\end{proof}
Our next step is to establish the bound of $\mathcal{J}^\prime(u_h, v_h)$ in the dual norm of $\mathbb{V}_{d}\times \mathbb{V}_c$, which is presented in the forthcoming lemma.
\begin{lemma}
Assume that all conditions of the Lemma \ref{4.5jbddlemma} holds. Then, for $(u_h, v_h)\in \mathbb{V}_{d}\times \mathbb{V}_c$, we have 
\begin{align}
\|\mathcal{J}^\prime(u_h, v_h)\|_{(\mathbb{V}_{d}\times \mathbb{V}_c)^\ast}&\leq ~ c_{46}\, \Big\{\sum_{\tau \in \mathscr{T}_h} |\eta_\tau(u_h,v_h)|^2 \Big\}^{\frac{1}{2}}, 
\end{align} 
where the constant $c_{46}$ is defined in Lemma \ref{4.5jbddlemma}.
\end{lemma}
\begin{proof}
We utilize the Cauchy-Schwarz inequality and the inequality \eqref{4.12lemmabdd} to derive 
\begin{align}
\|\mathcal{J}^\prime(u_h, v_h)\|_{(\mathbb{V}_{d}\times \mathbb{V}_c)^\ast}&=~ \sup_{(\psi, \varphi) \in \mathbb{V}_{d}\times \mathbb{V}_c} \frac{ ~|\mathcal{J}^\prime(u_h, v_h; \psi, \varphi)|~ }{ ~\big(\|\psi\|^2_{\mathbb{V}}+\|\varphi\|^2_{\mathbb{V}} \big)^{\frac{1}{2}}~}\nonumber\\
&\leq ~ c_{46}\, \Big\{\sum_{\tau \in \mathscr{T}_h} |\widetilde{\eta}_{\tau}(u_h,v_h)|^2+|\widehat{\eta}_{\tau}(u_h,v_h)|^2 \Big\}^{\frac{1}{2}}. 
\end{align} 
Setting $|\eta_\tau(u_h,v_h)|^2:= |\widetilde{\eta}_{\tau}(u_h,v_h)|^2+|\widehat{\eta}_{\tau}(u_h,v_h)|^2$, where the indicators $\widetilde{\eta}_{\tau}(u_h,v_h)$ and $\widehat{\eta}_{\tau}(u_h,v_h)$ are defined as in Lemma \ref{4.5jbddlemma}, this yields the desired result. 
\end{proof}
To facilitate the adaptive algorithm, we introduce a composite local indicator 
\begin{align}
\eta_\tau(u_h,v_h):= \big\{|\widetilde{\eta}_{\tau}(u_h,v_h)|^2+|\widehat{\eta}_{\tau}(u_h,v_h)|^2\big\}^{\frac{1}{2}}, \quad \text{for all} ~~ \tau \in \mathscr{T}_h,
\end{align}
where the indicators $\widetilde{\eta}_{\tau}$ and $\widehat{\eta}_{\tau}$ defined in equations \eqref{esteta4.11aa}-\eqref{esteta4.11bb}, which serve as the local refinement indicators for $u_h$ and $v_h$, respectively. 
In order to improve the accuracy of the numerical solution, this composite indicator will be a crucial part of our adaptive algorithm, directing the refining process.
\subsection{Post-optimization refinement technique} \label{algorithm1withtol} 
In order to start the adaptive algorithm, we define important parameters. Set tolerances $\Xi_v$, $\Xi_{v_n}$, and $\Xi_{RF}$ to stop the minimization and refinement loops, respectively; mesh size $h_k=\max_{\tau \in \mathscr{T}^n_{h_k}} diam(\tau)$ as the maximum diameter of elements in $n$-th time level and the $k$-th refinement $\mathscr{T}^n_{h_k}$; and marking parameter $0 < \vartheta \leq 1$ to determine elements for refinement. Then, using these parameters, we build an adaptive algorithm that uses the indicators to guide the refining mechanism and improve the accuracy of the solution.
\begin{flushleft}
\Hrule{4pt}{6pt} 
\vspace{0.1cm}
\textbf{Adaptive Algorithm-1} \vspace{0.1cm}
\Hrule{4pt}{6pt} \vspace{0.2cm}

\textbf{Step 1.~}{\tt Initialization:} Input crack field $v_0$~ and~ $\mathscr{T}_{h_0}$.  \vspace{0.2cm}

\textbf{Step 2.~}{\tt Start time loop:} For $t_n$, $n=1,\,2,\ldots,\, N_T.$ \\
\hspace{2cm} Compute $u^n_j$ and $v^n_j$ at $t_n$, for $j=1,\,2,\ldots$ \\
\hspace{2cm} Compute $u^n_j=\underset{\bar{u}\in \mathbb{V}_{f}}{{argmin}} \{\mathcal{J}(\bar{u}, v_j^{n-1})\}$ \\

\hspace{2cm} Compute $v_j^n=\underset{\bar{v}\in \mathbb{V}_{c}}{{argmin}} \{\mathcal{J}(u_j^n, \bar{v})\}$
\\ \vspace{0.2cm}

\hspace{2.0cm} \textbf{Check:} If $ \|v_j^{n}-v_j^{n-1}\|_{L^\infty(\mathcal{D})}\geq \Xi_{v_n}$\\

\hspace{5cm} Repeat \textbf{\tt Step 2} \\

\hspace{4.5cm} Else  \\

\hspace{5cm} \textbf{\tt Break}; \\

\hspace{4.5cm} End Else \\
\hspace{3.5cm}  If $v^n_j \leq \Xi_v $\\
\hspace{4.5cm} $v^n_j= 0$\\
\hspace{4.5cm} Else if $ v_n^j\geq 1.0$\\
\hspace{5cm} $v_n^j=1.0$\\
\hspace{4.5cm} End Else if\\
\hspace{4.0cm} End If \\
\hspace{3.5cm} End If \\ \vspace{0.2cm}
\textbf{Step 3.} Set $u_j=u_j^n$ and $v_j=v_j^n$. \\ \vspace{0.2cm}
\textbf{Step 4.} If $\big(\sum_{\tau \in \mathscr{T}^n_{h_j}}|\eta_{\tau}(u_j, v_j)|^2\big)^{1/2}> \Xi_{RF}$, \\ \vspace{0.1cm} 

\hspace{2cm} Determine a smallest subset $\mathcal{M}^n_{j}$ of $\mathscr{T}^n_{h_j}$ satisfying \\ \vspace{0.1cm} 
\hspace{2cm} $\sum_{\tau \in \mathcal{M}^n_{j}}|\eta_{\tau}(u_j, v_j)|^2 \geq \vartheta \sum_{\tau \in \mathscr{T}^n_{h_j}}|\eta_{\tau}(u_j, v_j)|^2$ \\ \vspace{0.2cm} 
\hspace{2cm} Refine the set  $\tau \in \mathcal{M}^n_{j}$, then generate  new mesh $\mathscr{T}^n_{h_{j+1}}$ (say)\\ \vspace{0.1cm} 
\hspace{1.55cm} End If $\big(\sum_{\tau \in \mathscr{T}^n_{h_j}}|\eta_{\tau}(u_j, v_j)|^2\big)^{1/2} \leq \Xi_{RF}$.\\ \vspace{0.2cm} 

\textbf{Step 6.} Set $u_h(t_n)=u_j$, $v_h(t_n)=v_j$, \\ \vspace{0.1cm}  \hspace{1.5cm} Repeat the steps. 
\Hrule{4pt}{6pt}
\end{flushleft} \vspace{0.2cm}
In the following lemma, we will state and establish the important features of the sequences acquired using \texttt{Adaptive Algorithm-1}.
\begin{lemma}\label{lemma4.6bddseq}
Let $\{(u_i,v_i)\}_{i=1}^\infty$ be a sequence generated via the {\tt Adaptive Algorithm-1} (Sec. \ref{algorithm1withtol}) such that $\{(u_i,v_i)\}_{i=1}^\infty \subseteq \mathbb{V}_{f,h_i}\times \mathbb{V}_{c,h_i}$, then the sequence $\{(u_i,v_i)\}_{i=1}^\infty$ holds the following properties
\begin{enumerate}
    \item[\tt (i)]{\tt Pointwise boundedness:} $0 \leq v_i(x) \leq 1$ on $\mathcal{D}$ for all $i \in \mathbb{N}$,
    \item[\tt (ii)]{\tt Sequence boundedness:} the sequence $\{(u_i,v_i)\}_{i=1}^\infty$ is bounded in $\mathbb{V}\times \mathbb{V}$.
\end{enumerate}   
\end{lemma}
\begin{proof}
Property (i) is a direct consequence of Proposition \ref{prop4.2:vh}. To prove property (ii), we substitute $u_i:=u_i^n$ and $v_i:=v_i^{n+1}$ into \eqref{ddofJAfem} and \eqref{ddofJbfem}, respectively, which results in
\begin{align}
 \mathcal{A}_{h_i}(v_i^n;u_i^n, \psi_i)=~0 \quad \psi_i \in \mathbb{V}_{d,h_i}, \quad \text{and} \quad  \mathcal{B}_{h_i}(u_i^n;v_i^{n+1}, \varphi_i)=~0 \quad \varphi_i \in \mathbb{V}_{c,h_i}.\label{4.26ahk}
\end{align}
From equation \eqref{4.26ahk}, we know that $\mathcal{A}_{h_i}(v_i^n;u_i^n, \psi_i)=0$. Setting $\psi_i=u_i^n-\pi_{h_i}f \in \mathbb{V}_{d,h_i}$, we then get utilizing the operator property of $\pi_h$,
\begin{align}
&~~\int_{\mathcal{D}}\dfrac{\big((1-\kappa)\pi_h((v_i^n)^2)+\kappa \big)}{\left[1 + \beta^{\alpha}\,\big((1-\kappa)\pi_h((v_i^n)^2)+\kappa \big)^\alpha \, |\nabla u_i^n|^{2\alpha}\right]^{\frac{1}{\alpha}+1}}\; \nabla u_i^n \cdot \nabla(u_i^n- \pi_{h_i}f)\, dx=~0,\nonumber\\
&\implies ~~~ \kappa\,\|\nabla u_i^n\|^2_{L^2(\mathcal{D})}
\,\leq~ \|\nabla u_i^n\|_{L^2(\mathcal{D})}\,\|\nabla \pi_{h_i} f\|_{L^2(\mathcal{D})}, \nonumber
\end{align}
and hence,
\begin{align}
&\|\nabla u_i^n\|_{L^2(\mathcal{D})}\,\leq~ \frac{1}{\kappa}\; \|\nabla f\|_{L^2(\mathcal{D})}. \nonumber
\end{align}
Since $0<\kappa \ll 1$ is a fixed number. This implies that the sequence $\{\|\nabla u_i\|_{L^2(\mathcal{D})}\}_{i=1}^\infty$ is bounded. 

Our objective is to demonstrate that the sequence $\{\|u_i\|_{L^2(\mathcal{D})}\}_{i=1}^\infty$ is bounded in $\mathbb{V}$. To establish this, we use the Friedrichs' inequality with constant $C_F>0$, which has
\begin{align}
\|u_i\|_{L^2(\mathcal{D})}&\leq~\|u_i-f\|_{L^2(\mathcal{D})}+\|f\|_{L^2(\mathcal{D})}\nonumber\\
&\leq~c_F\;\|\nabla(u_i-f)\|_{L^2(\mathcal{D})}+\|f\|_{L^2(\mathcal{D})}\nonumber\\
&\leq~c_F\; \|\nabla u_i\|_{L^2(\mathcal{D})}+(c_F^2+1)^{1/2}\,\|f\|_{\mathbb{V}}.\nonumber
\end{align}
Since $u_i-\pi_{h_i}f \in \mathbb{V}_d$ and $\pi_{h_i}f=f$.  This shows that $\{ u_i \}_{i=1}^\infty$ is a bounded sequence in $\mathbb{V}$. 

We show that the sequence $\{ v_i \}_{i=1}^\infty$ is a bounded sequence in $\mathbb{V}$. We go forward by setting $\varphi_i=v_i^{n+1} \in \mathbb{V}_{c,h_i}$ in equation \eqref{4.26ahk}, $\mathcal{B}_{h_i}(u_i^n;v_i^{n+1}, \varphi_i)=~0$, to acquire
\begin{align}
& \int_{\mathcal{D}}2\,\rho\, |\nabla v^{n+1}_{i} |^2\, dx+ \underaccent{ {\bf \geq 0}}{\underbrace{\int_{\mathcal{D}}\dfrac{(1-\kappa)\,|\nabla u^n_{i}|^2}{\left( 1 + \beta^{\alpha}\,\big((1-\kappa)\pi_{h_i}(v^{n+1}_{i})+\kappa \big)^\alpha\, |\nabla u^n_{i}|^{2\alpha}\right)^{\frac{1}{\alpha}+1}}\;\pi_h((v^{n+1}_{i})^2)\, dx}} \nonumber \\&~~~=~\int_{\mathcal{D}} \delta \,\pi_{h_i}(v^{n+1}_{i})\, dx. \nonumber 
\end{align}
Using the A.M-G.M inequality, we get
\begin{align}
&\int_{\mathcal{D}} |\nabla v^{n+1}_{i} |^2\, dx \leq~ \frac{\delta}{\rho}\int_{\mathcal{D}} v^{n+1}_{i} \,dx.\nonumber
\end{align}
Then, use of Cauchy-Schwarz inequality leads to
\begin{align}
\|\nabla v^{n+1}_{i} \|_{L^2(\mathcal{D})}\, dx \leq~ \frac{\delta\, \meas(\mathcal{D})}{\rho}.
\end{align}
This indicates that the sequence $\{\|\nabla v_i\|_{L^2(\mathcal{D})}\}_{i=1}^\infty$ is bounded. The sequence $\{ v_i \}_{i=1}^\infty$ is pointwise bounded, implying boundedness in the $L^2$-norm. Thus, the sequence $\{ v_i \}_{i=1}^\infty$ is bounded in $\mathbb{V}$. This concludes that the sequence $\{(u_i, v_i)\}_{i=1}^\infty$ is bounded in $\mathbb{V}\times \mathbb{V}$, which accomplished the proof.
\end{proof}
The above components lay the groundwork for our upcoming convergence analysis. The following theorem proves that, unaffected by the initial conditions chosen, \texttt{Adaptive Algorithm-1} will produce discrete solutions that approach a critical point of $\mathcal{J}(\cdot, \cdot)$
with a decreasing smaller tolerance $\Xi_{RF}$.
\begin{theorem}[Convergence upto a Tolerance] \label{Cngwithtol}
Let $\mathcal{D}\subset \mathbb{R}^d$ be an open bounded domain. Further, we assume that there exists a sequence $\{(u_i,v_i)\}_{i=1}^{\infty}$ in $\mathbb{V}_f\times \mathbb{V}_c$ with $v_i(x)\in [0, 1]$ for $a.e.\; x\in \mathcal{D}$, and for some $\gamma^1_i$ and $\gamma^2_i$ with $\gamma^1_i,\,\gamma^2_i \rightarrow 0$ as $i \rightarrow \infty$, such that
\begin{align}
&\mathcal{A}(v_i; u_i, \psi) \leq~ \gamma^1_i\; \|\nabla \psi\|_{L^2(\mathcal{D})}  \quad \;    \quad \text{and} \quad 
\mathcal{B}(u_i; v_i, \varphi) \leq~ \gamma^2_i\; \|\nabla \varphi\|_{L^2(\mathcal{D})}, \label{4.30abinq}
\end{align}
for all $\psi\in \mathbb{V}_d$ and $\varphi \in \mathbb{V}_c^\infty$. Again, we assume that the sequence $\{(u_i,v_i)\}_{i=1}^{\infty}$ is a bounded sequence in $\mathbb{V}\times \mathbb{V}$. Then, there exists a subsequence of $\{(u_i,v_i)\}_{i=1}^{\infty}$ and $(u,v)$ in $\mathbb{V}_f\times \mathbb{V}_c$ with $v(x)\in [0, 1]$  $a.e.\; x\in \mathcal{D}$ such that $u_i$ and $v_i$ converges strongly to $u$ and $v$ as $i \rightarrow \infty$, respectively, in $\mathbb{V}$. Additionally, $u$ and $v$ satisfy 
\begin{align}\label{4.6criticalpts}  
&\mathcal{A}(v; u, \psi) =~0  \quad \; \forall \,   \psi\in \mathbb{V}_d, \quad
\text{and} \quad \mathcal{B}(u; v, \varphi)=~0  \quad \; \forall\,   \varphi \in \mathbb{V}_c^\infty. 
\end{align}
Hence, the function $\mathcal{J}(\cdot,\cdot)$ has a critical point $(u,v)$ in $\mathbb{V}_f\times \mathbb{V}_c^\infty$.
\end{theorem}
\begin{proof} We will demonstrate this using a two-step process, with more sub-steps added as necessary.\smallskip

\noindent
\textbf{\tt Step 1.} First, we prove that there is a convergent subsequence of $\{(u_i,v_i)\}_{i=1}^{\infty}$ in $\mathbb{V}_f\times \mathbb{V}_c$. 
Note that the sequence $\{(u_i,v_i)\}_{i=1}^{\infty}$ is bounded in $\mathbb{V}\times \mathbb{V}$ according to Lemma \ref{lemma4.6bddseq}. Its boundedness demonstrates that a weakly convergent subsequence exists, because $\mathbb{V}$ is a Hilbert space. In particular, a subsequence (not relabelled) exists such that $(u_i,v_i) \overset{w}{\longrightarrow} (u,v)$ as $i\rightarrow \infty$ in $\mathbb{V}\times \mathbb{V}$. The fact that $\mathbb{V}_f$ is a closed and convex subset of $\mathbb{V}$ makes it noteworthy that it is also weakly closed. The weak limit $u$ is therefore an element of $\mathbb{V}$. We now
define a set  $$\mathbb{W}:=\{w\in \mathbb{V}_c :~  0\leq w(x)\leq 1~~ a.e. ~~ x\in \mathcal{D}\},$$
which is clearly a closed convex subset of $\mathbb{V}$. We observe that $0\leq v(x)\leq 1,\; a.e., \, x\in \mathcal{D}$, since $v_i\in \mathbb{W},\, \forall\, i\in \mathcal{N}$. The compact embedding of $H^1(\mathcal{D})$ in $ L^2(\mathcal{D})$ guarantees the strong convergence of the sequence $\{(u_{i},v_{i})\}_{i=1}^{\infty}$ to $(u,v)$ in $L^2(\mathcal{D})\times L^2(\mathcal{D})$.

In the following step, we will show that $(u,v)$ satisfies  $\mathcal{A}(v; u, \psi) =~0,\forall \,   \psi\in \mathbb{V}_d\; \text{and} \;\;\; \mathcal{B}(u; v, \varphi)=~0,$ $  \forall\,   \varphi \in \mathbb{V}_c^\infty$, respectively. \smallskip

\noindent
\texttt{Step 2. (i)} We begin by establishing that $\mathcal{A}(v; u, \psi)=0, \forall \, \psi\in \mathbb{V}_d.$  To do this, let us fix an arbitrary $\psi\in \mathbb{V}_d$ in \eqref{ddofJA}, we obtain
\begin{align}
\mathcal{A}(v; u, \psi)&=~ \int_{\mathcal{D}}\Big[\dfrac{\big((1-\kappa)\,v^2+\kappa \big)}{\left[1 + \beta^{\alpha}\,\big((1-\kappa)\,v^2+\kappa \big)^\alpha \, |\nabla u|^{2\alpha}\right]^{\frac{1}{\alpha}+1}}\; \nabla u \cdot \nabla \psi \Big] \, dx \nonumber \\ 
&=\int_{\mathcal{D}} \dfrac{\big((1-\kappa)\,v_i^2+\kappa)}{\left[1 + \beta^{\alpha}\,\big((1-\kappa)\,v^2+\kappa \big)^\alpha \, |\nabla u|^{2\alpha}\right]^{\frac{1}{\alpha}+1}}  \; \nabla u_i \cdot \nabla \psi\, dx\nonumber\\
&~~~+ \int_{\mathcal{D}} \dfrac{\big((1-\kappa)\,v^2+\kappa \big) }{\left[1 + \beta^{\alpha}\,\big((1-\kappa)\,v^2+\kappa \big)^\alpha \, |\nabla u|^{2\alpha}\right]^{\frac{1}{\alpha}+1}}\; \nabla (u-u_i) \cdot \nabla \psi \, dx \nonumber\\
&~~~+\int_{\mathcal{D}} \dfrac{(1-\kappa)\,(v^2-v_i^2)}{\left[1 + \beta^{\alpha}\,\big((1-\kappa)\,v^2+\kappa \big)^\alpha \, |\nabla u|^{2\alpha}\right]^{\frac{1}{\alpha}+1}}  \; \nabla u_i \cdot \nabla \psi\, dx \nonumber\\
&:=~\mathcal{X}_i+\mathcal{Y}_i+\mathcal{Z}_i. \label{xyzbb}
\end{align}   
It is necessary to demonstrate that each of the sequences $\mathcal{X}_i,\, \mathcal{Y}_i,\, \mathcal{Z}_i$ converges to zero as $i$ approaches infinity. \smallskip

\noindent
\textbf{(a) \tt Convergence of sequence $\mathcal{X}_i$.} By applying the inequality \eqref{4.30abinq} to the sequence $\mathcal{X}_i$, we may arrive at
\begin{align}
|\mathcal{X}_i|\leq \int_{\mathcal{D}} \dfrac{\big((1-\kappa)\,v_i^2+\kappa)}{\left[1 + \beta^{\alpha}\,\big((1-\kappa)\,v^2+\kappa \big)^\alpha \, |\nabla u|^{2\alpha}\right]^{\frac{1}{\alpha}+1}}  \; \nabla u_i \cdot \nabla \psi\, dx \leq~ \gamma^1_i\, \|\nabla \psi\|_{L^2(\mathcal{D})}.
\end{align}
Since $\gamma^1_i \rightarrow 0$ as $i \rightarrow \infty$, thus we have  $\mathcal{X}_i \rightarrow 0$ as $i\rightarrow\infty$. \smallskip

\noindent
\textbf{(b) \tt Convergence of sequence $\mathcal{Y}_i$.} Next, we investigate the convergence for the sequence $\mathcal{Y}_i$. Applying the weak convergence criteria to the sequence $\{\nabla u_i\}_{i=1}^\infty$, which converges to $\nabla u$ in $(L^2(\mathcal{D}))^d$, since $\big((1-\kappa)\,v^2+\kappa \big)\; \nabla \psi \in (L^2(\mathcal{D}))^d$. This implies that $\mathcal{Y}_i$ goes to zero as $i$ tends infinity.

 Next, we will examine convergence of the sequence $\mathcal{Z}_i$. \smallskip

\noindent
\textbf{(c) \tt Convergence of sequence $\mathcal{Z}_i$.} Invoking the Cauchy-Schwarz inequality, we examine the sequence $\mathcal{Z}_i$, which yields
\begin{align}
\big|\mathcal{Z}_i\big| &\leq~ \int_{\mathcal{D}} \dfrac{(1-\kappa)\,|v^2-v_i^2|}{\left[1 + \beta^{\alpha}\,\big((1-\kappa)\,v^2+\kappa \big)^\alpha \, |\nabla u|^{2\alpha}\right]^{\frac{1}{\alpha}+1}}  \; |\nabla u_i|\,|\nabla \psi|\, dx  \nonumber\\
&\leq~ 2\,(1-\kappa)\,\int_{\mathcal{D}} |v-v_i|\; |\nabla u_i|\,|\nabla \psi|\, dx.  \nonumber\\
&\leq~ 2\,(1-\kappa)\,\Big(\int_{\mathcal{D}} |v-v_i|\; |\nabla \psi|^2\,dx\Big)^{\frac{1}{2}} \times \Big(\int_{\mathcal{D}}|\nabla  u_i|^2\,dx\Big)^{\frac{1}{2}}.  \label{zeqbdd}
\end{align}
Since $|v-v_i|\leq 1$ and $1/\left[1 + \beta^{\alpha}\,\big((1-\kappa)\,v^2+\kappa \big)^\alpha \, |\nabla u|^{2\alpha}\right]^{\frac{1}{\alpha}+1}\leq 1$. To conclude, it is suffices to show that the term $\int_{\mathcal{D}} |v-v_i|\; |\nabla \psi|^2\,dx$ as $i \rightarrow \infty$ vanishes as $i$ tends to infinity.

To proceed, we assume that the sequence $\{v_i \}_{i=1}^\infty$ has a convergent subsequence  $ \{ v_{i_n} \}_{n=1}^\infty$, such that $v_{i_n} \longrightarrow v$ almost everywhere in $\mathcal{D}$. Furthermore, we assume that this subsequence satisfies, for fixed $\psi\in \mathbb{V}$,
$$\lim_{n \rightarrow \infty} \int_{\mathcal{D}} |v-v_{i_n}|\;|\nabla \psi|^2\,dx=\limsup_{i \rightarrow \infty} \int_{\mathcal{D}} |v-v_i|\;|\nabla \psi|^2\,dx.$$
Application of the Dominated Convergence theorem leads to
$$\lim_{n \rightarrow \infty} \int_{\mathcal{D}} |v-v_{i_n}|\;|\nabla \psi|^2\,dx=0.$$
Thus
$$\limsup_{i \rightarrow \infty} \int_{\mathcal{D}} |v-v_i|\;|\nabla \psi|^2\,dx=0,$$
and hence 
$$\int_{\mathcal{D}} |v-v_i|\;|\nabla \psi|^2\,dx=0,\quad \text{i} \rightarrow \infty.$$ 
Therefore, equation \eqref{zeqbdd} implies $|\mathcal{Z}_i|\longrightarrow 0$ as $i \rightarrow \infty$.
Last but not least, integrating these findings into equation \eqref{xyzbb} leads to $\mathcal{A}(v; u, \psi) =~0, \forall\, \psi\in \mathbb{V}_d$.

Next, we show that $\mathcal{B}(u; v, \varphi)=0$ is satisfied for all test functions $\varphi\in \mathbb{V}^\infty_c$.\smallskip

\noindent 
\textbf{(ii) \tt To show $\mathcal{B}(u; v, \varphi) =0, \forall \,   \varphi\in \mathbb{V}^\infty_c$.}  From \eqref{ddofJb}, we have
\begin{align}
\mathcal{B}(u; v, \varphi)&=~ \int_{\mathcal{D}}\Big[2\,\rho\, \nabla v \cdot \nabla \varphi - \delta\, \varphi + \dfrac{(1-\kappa)\,|\nabla u|^2}{\left( 1 + \beta^{\alpha}\,\big((1-\kappa)v^2+\kappa\big)^\alpha\, |\nabla u|^{2\alpha}\right)^{\frac{1}{\alpha}+1}}\;v \, \varphi \Big] \, dx\nonumber\\
&=~\int_{\mathcal{D}}\Big[2\,\rho\, \nabla (v-v_i) \cdot \nabla \varphi+ \Big(\dfrac{(1-\kappa)\,|\nabla u|^2}{\left[ 1 + \beta^{\alpha}\,\big((1-\kappa)v^2+\kappa\big)^\alpha\, |\nabla u|^{2\alpha}\right]^{\frac{1}{\alpha}+1}}\Big)\,(v-v_i)\,\varphi \Big]\,dx\nonumber\\
&~~~~+ \int_{\mathcal{D}}\Big[2\,\rho\, \nabla v_i \cdot \nabla \varphi -  \delta \, \varphi + \Big(\dfrac{(1-\kappa)\,|\nabla u_i|^2}{\left[ 1 + \beta^{\alpha}\,\big((1-\kappa)v^2+\kappa\big)^\alpha\, |\nabla u|^{2\alpha}\right]^{\frac{1}{\alpha}+1}}\Big)\,v_i\,\varphi \Big]\,dx\nonumber\\
&~~~~+\int_{\mathcal{D}} \dfrac{(1-\kappa)\,\big(|\nabla u|^2-|\nabla u_i|^2\big)}{\left[ 1 + \beta^{\alpha}\,\big((1-\kappa)v^2+\kappa\big)^\alpha\, |\nabla u|^{2\alpha}\right]^{\frac{1}{\alpha}+1}}\,v_i\,\varphi \, dx\nonumber\\
&:=~ \widetilde{\mathcal{X}}_i+\widetilde{\mathcal{Y}}_i+\widetilde{\mathcal{Z}}_i. \label{barxyzbb} 
\end{align}
We now show that each of the sequence  $\widetilde{\mathcal{X}}_i,\, \widetilde{\mathcal{Y}}_i,\, \widetilde{\mathcal{Z}}_i \longrightarrow 0$ as $i \rightarrow \infty$. We will examine each of them separately.\smallskip

\noindent
\textbf{(d) \tt Convergence of sequence $\widetilde{\mathcal{X}}_i$.} 
We start by analyzing the first term $\widetilde{\mathcal{X}}_i$, which yields
\begin{align}
|\widetilde{\mathcal{X}}_i|&= \Big|\int_{\mathcal{D}}\Big[2\,\rho\, \nabla (v-v_i) \cdot \nabla \varphi+ \Big(\dfrac{(1-\kappa)\,|\nabla u|^2}{\left[ 1 + \beta^{\alpha}\,\big((1-\kappa)v^2+\kappa\big)^\alpha\, |\nabla u|^{2\alpha}\right]^{\frac{1}{\alpha}+1}}\Big)\,(v-v_i)\,\varphi \Big]\,dx\Big| \nonumber\\
&\leq  2\,\rho\, \Big|\int_{\mathcal{D}} \nabla (v-v_i) \cdot \nabla \varphi\, dx \Big|+ \Big| \int_{\mathcal{D}} \big( (1-\kappa)\,|\nabla u|^2 \big)\, (v-v_i)\,\varphi \, dx \Big| \nonumber\\
&\leq  2\,\rho\, \Big|\int_{\mathcal{D}} \nabla (v-v_i) \cdot \nabla \varphi\, dx \Big|+ |1-\kappa|\,\|\varphi\|_{L^\infty(\mathcal{D})}\,\Big| \int_{\mathcal{D}} |\nabla u|^2\, (v-v_i)\, dx \Big|. 
\end{align}
As $1/ \big[ 1 + \beta^{\alpha}\,\big((1-\kappa)v^2+\kappa\big)^\alpha\, |\nabla u|^{2\alpha}\big]^{\frac{1}{\alpha}+1}\leq 1$. It is well acknowledged that $\nabla v_i \overset{w}{\longrightarrow} \nabla v$ and $v_i \longrightarrow  v$ in $(L^2(\mathcal{D}))^d$.  By an easy application of the preceding argumentation employed for ${\tt \mathcal{Z}_i}$ with $\psi=u$, it clearly follows that $\widetilde{\mathcal{X}}_i \rightarrow 0$ as $i \rightarrow \infty$.\smallskip

\noindent
\textbf{(e) \tt Convergence of sequence $\widetilde{\mathcal{Y}}_i$.} 
By using inequality \eqref{4.30abinq}, one can easily find that
$$|\widetilde{\mathcal{Y}}_i|\leq~ \gamma_i^2\,\|\varphi\|_{\mathbb{V}}.$$ 
Since $\gamma_i^2 \rightarrow 0$ as $ i \rightarrow \infty$, thus we obtain  $\widetilde{\mathcal{Y}}_i \longrightarrow 0$ as $ i \rightarrow \infty$.
Our focus now shifts to demonstrating the convergence of $\widetilde{\mathcal{Z}}_i$. \smallskip

\noindent
\textbf{(f) \tt Convergence of sequence $\widetilde{\mathcal{Z}}_i$.}  In order to do this, we first demonstrate that the sequence $\{\nabla u_i\}_{i=1}^\infty$ in $(L^2(\mathcal{D}))^d$ converges strongly, that is, $\nabla u_i \longrightarrow \nabla u$ as $i \rightarrow \infty$.\smallskip

\noindent 
\textbf{\tt {Strong convergence for the sequence $\{\nabla u_i\}_{i=1}^\infty$ in $(L^2(\mathcal{D}))^d$}.} Taking into account that $u-u_i\in \mathbb{V}_d$, we arrive at
\begin{align}
\kappa \, \|\nabla u- \nabla u_i\|^2_{L^2(\mathcal{D})}&\leq~ \int_{\mathcal{D}} \dfrac{\big((1-\kappa) v_i^2+\kappa \big)}{\left[ 1 + \beta^{\alpha}\,\big((1-\kappa)v^2+\kappa\big)^\alpha\, |\nabla u|^{2\alpha}\right]^{\frac{1}{\alpha}+1}}\,\,(\nabla u- \nabla u_i)\cdot (\nabla u- \nabla u_i)\, dx \nonumber\\
&=~\int_{\mathcal{D}} \dfrac{\big((1-\kappa) v_i^2+\kappa \big)}{\left[ 1 + \beta^{\alpha}\,\big((1-\kappa)v^2+\kappa\big)^\alpha\, |\nabla u|^{2\alpha}\right]^{\frac{1}{\alpha}+1}}\,(- \nabla u_i)\cdot (\nabla u- \nabla u_i)\, dx \nonumber\\
&~~~~+ \int_{\mathcal{D}} \dfrac{\big((1-\kappa) v_i^2+\kappa \big)}{\left[ 1 + \beta^{\alpha}\,\big((1-\kappa)v^2+\kappa\big)^\alpha\, |\nabla u|^{2\alpha}\right]^{\frac{1}{\alpha}+1}}\,\nabla u\cdot (\nabla u- \nabla u_i)\, dx.\nonumber
\end{align}
Utilizing the inequality \eqref{4.30abinq}, the condition $\mathcal{A}(v;u,\psi)=0$ with $\psi=u-u_i$, and the Cauchy-Schwarz inequality, we derive 
\begin{align}
\kappa \, \|\nabla u- \nabla u_i\|^2_{L^2(\mathcal{D})}&\leq~ \gamma^1_i\,\|\nabla (u-u_i)\|_{L^2(\mathcal{D})}+ \int_{\mathcal{D}} \dfrac{(1-\kappa)\, (v_i^2-v^2)\,\nabla u\cdot (\nabla u- \nabla u_i)}{\left[ 1 + \beta^{\alpha}\,\big((1-\kappa)\,v^2+\kappa\big)^\alpha\, |\nabla u|^{2\alpha}\right]^{\frac{1}{\alpha}+1}}\, dx\nonumber\\
&\leq~ \gamma^1_i\,\|\nabla (u-u_i)\|_{L^2(\mathcal{D})}+ c_s\Big(\int_{\mathcal{D}} 4\,|1-\kappa||v-v_i|\,|\nabla u|^2 dx\Big)^{\frac{1}{2}} \nonumber\\
& ~~~~ \times \Big(\int_{\mathcal{D}} |\nabla u- \nabla u_i|^2\,dx\Big)^{\frac{1}{2}}.\nonumber
\end{align}
Notice that $0\leq v,\,v_i\leq 1,~ a.e.,~ x\in \mathcal{D}$, $|v-v_i|\leq 1$, $1/\left[ 1 + \beta^{\alpha}\,\big((1-\kappa)\,v^2+\kappa\big)^\alpha\, \|\nabla u\|^{2\alpha}\right]^{\frac{1}{\alpha}+1}\leq 1,$  and $|v^2-v_i^2|^2=|v-v_i|^2|v+v_i|^2\leq 4\,|v-v_i|$,  where $c_s$ is the Cauchy-Schwarz inequality constant, and hence 
\begin{align}
\kappa \, \|\nabla u- \nabla u_i\|_{L^2(\mathcal{D})}&\leq~ \gamma^1_i+ 2\,c_s\, |1-\kappa|\,\Big(\int_{\mathcal{D}}|v-v_i|\,|\nabla u|^2 dx\Big)^{\frac{1}{2}}. \label{4.37gradu}
\end{align}
This leads to $\lim_{i \rightarrow \infty}\|\nabla u- \nabla u_i\|_{L^2(\mathcal{D})}=0$, because  $\gamma^1_i\longrightarrow 0$ as $i \rightarrow \infty$ and  $\int_{\mathcal{D}}|v-v_i|\,|\nabla u|^2 dx \longrightarrow 0$ from the previous step of ${\tt \mathcal{Z}_i}$ with $\psi=u$.

Let us continue our investigation to the convergence of $\widetilde{\mathcal{Z}}_i$, commencing with
\begin{align}
|\widetilde{\mathcal{Z}}_i|&\leq~ \int_{\mathcal{D}} \Big| \dfrac{(1-\kappa)\,\big(|\nabla u|^2-|\nabla u_i|^2\big)}{\left[ 1 + \beta^{\alpha}\,\big((1-\kappa)v^2+\kappa\big)^\alpha\, |\nabla u|^{2\alpha}\right]^{\frac{1}{\alpha}+1}}\,v_i\,\varphi \Big| \, dx\nonumber\\
& \leq |1-\kappa|\,\|\varphi\|_{L^\infty(\mathcal{D})}\int_{\mathcal{D}} \dfrac{\big| |\nabla u|^2-|\nabla u_i|^2\big|\;|v_i|}{\big| \left[ 1 + \beta^{\alpha}\,\big((1-\kappa)v^2+\kappa\big)^\alpha\, |\nabla u|^{2\alpha}\right]^{\frac{1}{\alpha}+1}\big|} \, dx\nonumber\\
& \leq |1-\kappa|\,\|\varphi\|_{L^\infty(\mathcal{D})}\|\nabla u-\nabla u_i\|_{L^2(\mathcal{D})}\,\|\nabla u+\nabla u_i\|_{L^2(\mathcal{D})}.\nonumber 
\end{align}
Since $|v_i|\leq 1$ and $1/\big| \left[ 1 + \beta^{\alpha}\,\big((1-\kappa)v^2+\kappa\big)^\alpha\, \|\nabla u\|^{2\alpha}\right]^{\frac{1}{\alpha}+1}\big| \leq 1$. The convergence of  $\{\nabla u_i\}_{i=1}^\infty$ to $\nabla u$ in $(L^2(\mathcal{D}))^d$ implies that  $\widetilde{\mathcal{Z}}_i \longrightarrow 0$ as $i\rightarrow \infty$. 
We finally conclude that $\mathcal{B}(u; v, \varphi)=0,$ $ \forall \, \varphi\in \mathbb{V}^\infty_c$ by combining the convergence results $\widetilde{\mathcal{X}}_i,\, \widetilde{\mathcal{Y}}_i$ and $\widetilde{\mathcal{Z}}_i$, respectively.

The last step is to prove that the sequence $\{\nabla v_i\}_{i=1}^\infty$ in $(L^2(\mathcal{D}))^d$ converges strongly. \smallskip

\noindent 
{\tt {Strong convergence for the sequence $\{\nabla v_i\}_{i=1}^\infty$ in $(L^2(\mathcal{D}))^d$}.} To demonstrate the strong convergence of $\{\nabla v_i\}_{i=1}^\infty$, we take into
\begin{align}
&2\,\rho\, \|\nabla v_i-\nabla v\|^2_{L^2(\mathcal{D})}\leq~\int_{\mathcal{D}}  \Big(\dfrac{(1-\kappa)\,|\nabla u|^2}{\left[ 1 + \beta^{\alpha}\,\big((1-\kappa)v^2+\kappa\big)^\alpha\, |\nabla u|^{2\alpha}\right]^{\frac{1}{\alpha}+1}}\Big)\,|v_i-v|^2\,dx \nonumber\\
&~~~~~+\int_{\mathcal{D}} 2\,\rho\, (\nabla v_i-\nabla v)\cdot  (\nabla v_i-\nabla v)\,dx \nonumber\\
&~~= \int_{\mathcal{D}}  \Big[\dfrac{(1-\kappa)\,|\nabla u|^2}{\left[ 1 + \beta^{\alpha}\,\big((1-\kappa)v^2+\kappa\big)^\alpha\,|\nabla u|^{2\alpha}\right]^{\frac{1}{\alpha}+1}}\,v_i\,(v_i-v)+ 2\,\rho\,\nabla v_i\cdot  (\nabla v_i-\nabla v)\Big]\,dx\nonumber\\
&~~~~~ -\int_{\mathcal{D}}  \Big[\dfrac{(1-\kappa)\,|\nabla u|^2}{\left[ 1 + \beta^{\alpha}\,\big((1-\kappa)v^2+\kappa\big)^\alpha\,|\nabla u|^{2\alpha}\right]^{\frac{1}{\alpha}+1}}\,v\,(v_i-v)+ 2\,\rho\,\nabla v \cdot  (\nabla v_i-\nabla v)\Big]\,dx \nonumber\\
&~~\leq \mathcal{B}(u_i;v_i,v_i-v)-\mathcal{B}(u;v,v_i-v) +\int_{\mathcal{D}}  (1-\kappa)\,\big(|\nabla u|^2-|\nabla u_i|^2\big)\,v_i\,(v_i-v)\, dx \nonumber\\
&~~ \leq  \mathcal{B}(u_i;v_i,v_i-v)-\mathcal{B}(u;v,v_i-v) + |1-\kappa|\|v_i\|_{L^\infty(\mathcal{D})}\,\|\nabla u- \nabla u_i\|_{L^2(\mathcal{D})}\,\nonumber\\
&~~~~~~\times \|\nabla u+\nabla u_i\|_{L^2(\mathcal{D})}\;\|v_i-v\|_{L^2(\mathcal{D})}. \nonumber
\end{align}
By applying $\mathcal{B}(u;v,v_i-v)=0$, for $v_i-v \in \mathbb{V}_c^\infty$, together with \eqref{4.30abinq} and \eqref{4.37gradu}, we acquire
\begin{align}
2\,\rho\, \|\nabla v_i-\nabla v\|^2_{L^2(\mathcal{D})}
&\leq~  \gamma^2_i\,\|\nabla(v_i-v)\|_{L^2(\mathcal{D})}+ \Big|\frac{1-\kappa}{\kappa} \Big|\|v_i\|_{L^\infty(\mathcal{D})}\,\|\nabla u+\nabla u_i\|_{L^2(\mathcal{D})}\nonumber\\
&~~~~\times \|v_i-v\|_{L^2(\mathcal{D})} \, \Big(\gamma^1_i+ 2\,c_s\, |1-\kappa|\,\big(\int_{\mathcal{D}}|v-v_i|\,|\nabla u|^2 dx\big)^{\frac{1}{2}}\Big), \nonumber 
\end{align}
and hence
\begin{align}
\|\nabla v_i-\nabla v\|_{L^2(\mathcal{D})}&\leq~ \frac{1}{2\,\rho}\Big[\gamma^2_i+ \Big|\frac{1-\kappa}{\kappa} \Big|\|v_i\|_{L^\infty(\mathcal{D})} \Big(\gamma^1_i+ 2\,c_s\, |1-\kappa|\,\big(\int_{\mathcal{D}}|v-v_i|\,|\nabla u|^2 dx\big)^{\frac{1}{2}}\Big)\nonumber\\
&~~~~\times \|\nabla u+\nabla u_i\|_{L^2(\mathcal{D})} \Big]. \nonumber
\end{align}
Since $\gamma^1_i,\, \gamma^2_i \longrightarrow 0$ as $i \rightarrow \infty$, and $\int_{\mathcal{D}}|v_i-v|\,|\nabla u|^2 dx \longrightarrow 0$. These convergence results, along with the boundedness of $\{\|\nabla u+\nabla u_i\|_{L^2(\mathcal{D})}\}_{i=0}^\infty$ and $\{\|\nabla(v_i-v)\|_{L^2(\mathcal{D})}\}_{n=1}^\infty$, we determine that $\{\nabla v_i\}_{i=1}^\infty$  is strongly converges to $\nabla v$ in $(L^2(\mathcal{D}))^d$ as $i \rightarrow \infty$. This concludes the proof of the theorem.
\end{proof}
The following section addresses a modified adaptive approach and conducts a convergence analysis.
\subsection{Iterative optimization refinement technique}
We specified the following settings to govern the adaptive algorithm. We set the tolerances
{\tt $\Xi_v$ and $\Xi_{RF}$}, respectively, to stop the minimization procedure and the refinement loop in the following algorithm. Furthermore, {\tt $\mathscr{T}^n_{h_0}$} denotes the  $n$-th time level with initial mesh $\mathscr{T}_{h_0}$, however, 
{\tt $h^n_0$} refers to the initial mesh size at the $n$-th time level, which is a maximum diameter of elements in $\mathscr{T}^n_{h_0}$. In addition, we propose the marking parameter {\tt $\vartheta \in (0, 1]$}, which is used to identify the set of elements to be refined. Note that the mesh is refined in each subsequent minimization step ($n=k/2$, where $k\in \mathbb{N}$) in the adaptive refinement procedure. Further, $\mathscr{T}^n_{h_k}$ represents the $k$-th refinement of the $n$-th time level, with a mesh size of $h^n_k=\max_{\tau \in \mathscr{T}^n_{h_k}} diam(\tau)$. For simplicity, we write $h_k$ in place of $h_k^n$. The refining procedure employs the bisection technique to guarantee that the final mesh satisfies the shape-regularity constraints \cite{Mitchell1989}. 
Below is a description of the adaptive algorithm used to solve the problem:
\begin{flushleft}
\Hrule{6pt}{4pt} 
\vspace{0.1cm}
\textbf{Adaptive Algorithm-2} \vspace{0.1cm} \label{adaptalgo2}
\Hrule{6pt}{4pt} \vspace{0.1cm}

\textbf{Step 1.~}{\tt Initialization:} Input crack field $v_0$~ and~ $\mathscr{T}^{1/2}_{h_0}$.  \\
\vspace{0.25cm}

\textbf{Step 2.~} For $n=1,\,2,\ldots,\,.$\\  \vspace{0.1cm}
\hspace{2cm} Set $\mathscr{T}^n_{h_0}=\mathscr{T}^{n-1/2}_{h_0}$\\  \vspace{0.25cm}

\hspace{1cm} \textbf{A.~}{\tt Start time loop for $u$:} For $t_n$, $n=1,\,2,\ldots,\, N_T.$ \\ \vspace{0.1cm}
\hspace{2cm} Compute $u^n_k$ for $k=1,\,2,\ldots\,$ \\  \vspace{0.1cm}
\hspace{2cm} by using the formula $u^n_k=\underset{\bar{u}\in \mathbb{V}_{f, h^n_k}}{{argmin}} \{\mathcal{J}_{h^n_k}(\bar{u}, v_k^{n-1})\}$ \\ \vspace{0.25cm}

\hspace{1cm}\textbf{Check:} If $\big(\sum_{\tau \in \mathscr{T}^n_{h_k}}|\eta_{\tau}(u_k^n, v_k^{n-1})|^2\big)^{1/2}> \frac{\Xi_{RF}}{\sqrt{2}}$, \\ \vspace{0.1cm} 

\hspace{3cm} Determine a smallest subset $\mathcal{M}^n_k$ of $\mathscr{T}^n_{h_k}$ satisfying \\ \vspace{0.1cm} 
\hspace{3cm} $\sum_{\tau \in \mathcal{M}^n_k}|\eta_{\tau}(u_k^n, v_k^{n-1})|^2 \geq \vartheta \sum_{\tau \in \mathscr{T}^n_{h_k}}|\eta_{\tau}(u_k^n, v_k^{n-1})|^2$ \\ \vspace{0.2cm} 
\hspace{3cm} Refine the set $\tau \in \mathcal{M}^n_k$ and then generate new mesh $\mathscr{T}^n_{h_{k+1}}$ \\ \vspace{0.1cm} 
\hspace{3.55cm} End If $\big(\sum_{\tau \in \mathscr{T}^n_{h_k}}|\eta_{\tau}(u_k^n, v_k^{n-1})|^2\big)^{1/2} \leq \frac{\Xi_{RF}}{\sqrt{2}}$.\\ \vspace{0.2cm} 

\hspace{3cm} Set $u^n=u^n_k, \mathscr{T}^n_{k}= \mathscr{T}^n_{h_k}$ and $\mathscr{T}^{n+1/2}_{k,1}=\mathscr{T}^n_{k}$. \\ \vspace{0.3cm}

\hspace{1cm} \textbf{B.~}{\tt Start time loop for $v$:} For $t_n$, $n=1,\,2,\ldots,\, N_T.$\\
\hspace{2cm} Compute $v^n_k$ for $k=1,\,2,\ldots\,$ \\
\hspace{2cm} by using the formula $v_k^n=\underset{\bar{v}\in \mathbb{V}_{c,h_k^{n+1/2}}}{{argmin}} \{\mathcal{J}_{h_k^{n+1/2}}(u_k^n, \bar{v})\}$\\ \vspace{0.2cm}

\hspace{1cm}\textbf{Check:} If $\big(\sum_{\tau \in \mathscr{T}^n_{h_k}}|\eta_{\tau}(u_k^n, v_k^n)|^2\big)^{1/2}> \frac{\Xi_{RF}}{\sqrt{2}}$, \\ \vspace{0.1cm} 

\hspace{3cm} Determine a smallest subset $\mathcal{M}^{n+1/2}_k$ of $\mathscr{T}^{n+1/2}_{h_k}$ satisfying \\ \vspace{0.1cm} 
\hspace{3cm} $\sum_{\tau \in \mathcal{M}^{n+1/2}_k}|\eta_{\tau}(u_k^n, v_k^n)|^2 \geq \vartheta \sum_{\tau \in \mathscr{T}^{n+1/2}_{h_k}}|\eta_{\tau}(u_k^n, v_k^n)|^2$ \\ \vspace{0.2cm} 
\hspace{3cm} Refine the set  $\tau \in \mathcal{M}^{n+1/2}_{h_k}$ and then generate new mesh $\mathscr{T}^{n+1/2}_{h_{k+1}}$ \\ \vspace{0.1cm} 
\hspace{3.55cm} End If $\big(\sum_{\tau \in \mathscr{T}^{n+1/2}_{h_{k+1}}}|\eta_{\tau}(u_k^n, v_k^n)|^2\big)^{1/2} \leq \frac{\Xi_{RF}}{\sqrt{2}}$.\\ 
\vspace{0.25cm}
\hspace{3cm} Set $v^n=v^n_{h_k}$, and $\mathscr{T}^{n+1/2}_{k,1}=\mathscr{T}^{n+1}_{h_k}$. \\ \vspace{0.3cm}

\hspace{1.5cm} \textbf{Check:} If $ \|v^n-v^{n-1}\|_{L^\infty(\mathcal{D})}\geq \Xi_v$\\ \vspace{0.1cm}

\hspace{4.7cm} Repeat \textbf{\tt Step 2} \\\vspace{0.2cm}

\hspace{3cm} Else  \\ \vspace{0.1cm}

\hspace{3.7cm} \textbf{\tt Break}; \\ \vspace{0.1cm}

\hspace{3cm} End Else \\ \vspace{0.2cm}

\textbf{Step 4.} Set $u_{h_k}(t_n)=u^n$ and $v_{h_k}(t_n)=v^n$. \\ \vspace{0.1cm}  \hspace{2.5cm} Repeat the process. 
\Hrule{6pt}{3pt}
\end{flushleft} \vspace{0.1cm}
Within a certain tolerance, the series of discrete solutions produced in Step 2 of {\tt Adaptive Algorithm-2} converges to exact solutions. By supplying a more detailed mesh and permitting a lower regularization value, adaptive mesh refinement makes it possible to evaluate fracture paths with more accuracy. If the refinement tolerance is small enough, this local refinement capability is activated when a new crack is introduced in a coarse domain region. We will now establish crucial features in the improved {\tt Adaptive Algorithm-3}.\\

\noindent
\textbf{Improved Refinement Adaptive Algorithm-3:} \label{madpiii} Based on {\tt Algorithm-2} (Sec. \ref{adaptalgo2}), we propose Algorithm-3, which includes a revised refinement tolerance technique. At each $k$-th step, we substitute the constant refinement tolerance $\Xi_{(RF)}$ with a step-dependent tolerance $\Xi_{(RF)_k}$, which approaches zero as $k\rightarrow \infty$. The termination step in the {\tt Adaptive Algorithm-2} can be skipped if the condition $\|v^n-v^{n-1}\|_{L^\infty(\mathcal{D})}< \Xi_v$ is satisfied (see Sec. \ref{adaptalgo2}).

In the subsequent lemma, we develop and demonstrate the essential characteristics of the sequence generated by the improved {\tt Adaptive Algorithm-3}.
\begin{lemma}\label{410lmadpalg}
Let $\{(u^i,v^i)\}_{i=1}^\infty$ be a sequence generated by the {\it Adaptive Algorithm-3} (Sec. \ref{algorithm1withtol}) such that $\{(u^i,v^i)\}_{i=1}^\infty \subseteq \mathbb{V}_{f,h_i}\times \mathbb{V}_{c,h_i}$, then the sequence $\{(u^i,v^i)\}_{i=1}^\infty$ holds the following properties
\begin{enumerate}
 \item[\tt (i)]{\tt Pointwise boundedness of $v^i$:} $0 \leq v^i(x) \leq 1$ on $\mathcal{D}$ for all $i \in \mathbb{N}$,
\item[\tt (ii)]{\tt Sequence boundedness :}  $\{(u^i,v^i)\}_{i=1}^\infty$ is a bounded sequence in $\mathbb{V}\times \mathbb{V}$,
\item[\tt (iii)]{\tt Monotonicity condition:}  $\mathcal{J}_{h_i}(u^i, v^i)\leq \mathcal{J}_{h_{i-1/2}}(u^{i-1}, v^i)\leq \mathcal{J}_{h_{i-1}}(u^{i-1}, v^{i-1}),~ \forall i\in \mathbb{N}\; \text{with}\; i\geq 2,$
\item[\tt (iv)]{\tt Energy convergence I:}   $  \lim_{i\rightarrow \infty}|\mathcal{J}(u^i, v^i)-\mathcal{J}_{h_{i+1/2}}(u^i, v^i)|=0$,
\item[\tt (v)]{\tt Energy convergence II:}  $  \lim_{i\rightarrow \infty}|\mathcal{J}(u^i, v^i)-\mathcal{J}_{h_i}(u^i, v^i)|=0$,
\item[\tt (vi)]{\tt Energy descent condition:} $$\liminf_{i \rightarrow \infty}\mathcal{J}(u^i, v^i)\leq \liminf_{i \rightarrow \infty}\mathcal{J}(u^{i-1}, v^i)\leq \liminf_{i\rightarrow \infty}\mathcal{J}(u^{i-1}, v^{i-1}).$$
\end{enumerate}
\end{lemma}
\begin{proof}
The proofs of properties (i) and (ii) parallel those of the lemma \ref{lemma4.6bddseq}. We will now focus on establishing property (iii). \smallskip

\noindent 
\textbf{(iii) \tt Monotonicity condition.} The proof of property (iii) is based on the following essential observation:
\begin{align*}
u^i=\argmin_{\bar{u}\in \mathbb{V}_{f, h_i}}\{\mathcal{J}_{h_i}{(\bar{u}, v^i)}\}, 
\end{align*}
and hence
\begin{align}
\mathcal{J}_{h_i}{(u^i, v^i)}-\mathcal{J}_{h_i}{(u^{i-1}, v^i)}\leq~\frac{1}{2}\int_{\mathcal{D}} \big( ( 1- \kappa)\; \pi_h((v^i)^2) + \kappa  \big)\,\big(|\nabla u^i|^2 -|\nabla u^{i-1}|^2\big)\, dx. \label{43888jhijhi-1}
\end{align}
Observe that  $1/\big[1 + \beta^{\alpha}\,\big( ( 1- \kappa)\; \pi_{h_i}((v^i)^2) + \kappa  \big)^\alpha\, \|\nabla u^i\|^{2\alpha}\big]^{1/\alpha} \leq 1,\, \forall\,i$, and $(v^i)^2$ is a convex function which lies between $0$ and $1$. Further, $\pi_{h_i}((v^i)^2)$ is its piecewise linear interpolant. Additionally, the sequence $\{\|\nabla u^i\|\}_{i=1}^\infty$ is bounded, implying that
\begin{align}
 \mathcal{J}_{h_i}{(u^i, v^i)}-\mathcal{J}_{h_i}{(u^{i-1}, v^i)}\leq~0. \label{4.388ui}   
\end{align}
Furthermore, on each element $\tau\in \mathscr{T}_{h_{i-1/2}}$, the interpolation function of $(v^i)^2$ is $\pi_{h_i}((v^i)^2)$, which has a greater or equal degree of accuracy than $\pi_{h_{i-1/2}}((v^i)^2)$ owing to the finer mesh resolution. Which brings us to
\begin{align}
(v^i)^2 \leq~ \pi_{h_i}((v^i)^2) \leq~ \pi_{h_{i-1/2}}((v^i)^2). \label{4.38intv}
\end{align}
Invoking the bound provided by \eqref{4.38intv}, we apply similar approach as in \eqref{43888jhijhi-1} to determine that
\begin{align}
\mathcal{J}_{h_i}{(u^{i-1}, v^i)}-\mathcal{J}_{h_{i-1/2}}{(u^{i-1}, v^i)}\leq~0. \label{441jhijh}
\end{align}
Adding equations \eqref{4.388ui} and \eqref{441jhijh}, to obtain
\begin{align}
 \mathcal{J}_{h_i}{(u^i, v^i)}-\mathcal{J}_{h_{i-1/2}}{(u^{i-1}, v^i)}&=~\big( \mathcal{J}_{h_i}{(u^i, v^i)}-\mathcal{J}_{h_i}{(u^{i-1}, v^i)}\big)\nonumber\\ &~~~+\big(\mathcal{J}_{h_i}{(u^{i-1}, v^i)}-\mathcal{J}_{h_{i-1/2}}{(u^{i-1}, v^i)}\big)\leq~0, \nonumber 
\end{align}
and, hence $\mathcal{J}_{h_i}{(u^i, v^i)}\leq~\mathcal{J}_{h_{i-1/2}}{(u^{i-1}, v^i)}$. In a similar fashion, one can easily prove that $\mathcal{J}_{h_{i-1/2}}{(u^{i-1}, v^i)} \leq~ \mathcal{J}_{h_{i-1}}{(u^{i-1}, v^{i-1})},\, i\geq 2$. Hence, property (iii) is proven. \smallskip

\noindent
\textbf{(iv) \tt Energy convergence I.} Let us consider
\begin{align}
|\mathcal{J}{(u^i, v^{i+1})}-\mathcal{J}_{h_{i+1/2}}{(u^i, v^{i+1})}|&\leq ~\Big|\int_{\mathcal{D}}\Bigg(\dfrac{( 1- \kappa)|\nabla u^i|^2\, \big((v^{i+1})^2 -\pi_{h_{i+1/2}}(v^{i+1})^2\big)}{\big[ 1 + \beta^{\alpha}\,\big((1-\kappa)(v^{i+1})^2+\kappa \big)^\alpha\, | \nabla u^i|^{2\alpha}\big]^{\frac{1}{\alpha}+1}} \nonumber\\
&~~~~-\delta\,\big(v^{i+1} -\pi_{h_{i+1/2}}(v^{i+1}) \big) \Bigg)\, dx \Big|\nonumber\\
&~\leq~ \big|\mathcal{B}(u^i; v^{i+1}, v^{i+1})- \mathcal{B}_{h_{i+1/2}}(u^i; v^{i+1}, v^{i+1})\big|. \nonumber
\end{align}
Using the estimation method as $I_4$ in the equation \eqref{4.20BB} with $u_h=u^i$, $v_h=v^{i+1}$, and $\varphi_h=v^{i+1}$, the residual estimate $\varphi \mapsto \mathcal{B}(u^i; v^{i+1}, \varphi)$ with $\varphi=v^{i+1}$ produces the following
\begin{align}
|\mathcal{J}{(u^i, v^{i+1})}-\mathcal{J}_{h_{i+1/2}}{(u^i, v^{i+1})}|
&~\leq~ c_{49}\,\widehat{\widetilde{\eta}}_{h_{i+1/2}}(u^i,v^{i+1})\;\|\nabla v^{i+1}\|_{L^2(\mathcal{D})} \label{4.42jhi12}
\end{align}
where 
\begin{align}
\widehat{\widehat{\eta}}_{h_{i+1/2}}(u^i,v^{i+1})&=~\Bigg\{\sum_{ \tau_i\in \mathscr{T}_h} \Bigg(\int_{\tau_i} h_{\tau_i}^4\,\Big|\dfrac{(1-\kappa)\,|\nabla u^i|^2}{\left[ 1 + \beta^{\alpha}\,\big((1-\kappa)\,(v^{i+1})^2 +\kappa \big)^\alpha\, | \nabla u^i|^{2\alpha}\right]^{\frac{1}{\alpha}+1}} \Big|^2\,dx \Bigg)\nonumber \\
&~~~~~~\times \|\nabla v^{i+1}\|^2_{L^\infty(\tau_i)}  \Bigg\}^{\frac{1}{2}},\label{etahathat}
\end{align}
and the constant $c_{49}$ depends on $c_\varrho,\, c_{eq}$ and $c_{45}$, respectively. The implication is that $h_{i+1/2}$ tends to zero as $i$ approaches infinity while $\lim_{i \rightarrow \infty} \Xi_{(RF)_i}=0$. Consequently, as $i$ goes to $\infty$, the term $\widehat{\widehat{\eta}}_{h_{i+1/2}}(u^i,v^{i+1})$ vanishes. Therefore,
\begin{align*}
|\mathcal{J}{(u^i, v^{i+1})}-\mathcal{J}_{h_{i+1/2}}{(u^i, v^{i+1})}| \longrightarrow 0 \quad \text{as} \quad i \rightarrow \infty.
\end{align*}
This validates the property $(iv)$. We will next verify the property $(v)$.\smallskip

\noindent
\textbf{(v) \tt Energy convergence II.} We repeat the previous argument by using the condition $a.b\leq \frac{1}{2}(a^2+b^2)$ to determine
\begin{align}
&|\mathcal{J}{(u^i, v^i)}-\mathcal{J}_{h_i}{(u^i, v^i)}|\leq ~\Big|\int_{\mathcal{D}}\frac{1}{2}\, \bigg(\dfrac{( 1- \kappa)\,|\nabla u^i|^2\,\big((v^i)^2 -\pi_{h_i}(v^i)^2 \big)}{\big[ 1 + \beta^{\alpha}\,\big((1-\kappa)(v^i)^2+\kappa \big)^\alpha\, | \nabla u^i|^{2\alpha}\big]^{\frac{1}{\alpha}+1}}\nonumber\\
&\hspace{4.5cm} -2\,\delta \big(v^i -\pi_{h_i}(v^i)\big)\bigg)\,dx \Big|\nonumber\\
&\leq ~\Big|\int_{\mathcal{D}}\frac{1}{2}\, \dfrac{( 1- \kappa)\,\nabla u^i\cdot \nabla (u^i-u^{i-1})}{\big[ 1 + \beta^{\alpha}\,\big((1-\kappa)(v^i)^2+\kappa \big)^\alpha\, | \nabla u^i|^{2\alpha}\big]^{\frac{1}{\alpha}+1}} \times \big((v^i)^2 -\pi_{h_i}(v^i)^2 \big)\, dx \nonumber\\
&~~~~+\int_{\mathcal{D}}\frac{1}{2}\,\dfrac{( 1- \kappa)\,\nabla u^i\cdot \nabla u^{i-1} }{\big[ 1 + \beta^{\alpha}\,\big((1-\kappa)(v^i)^2+\kappa \big)^\alpha\, | \nabla u^i|^{2\alpha}\big]^{\frac{1}{\alpha}+1}} \times \big((v^i)^2 -\pi_{h_i}(v^i)^2 \big) \, dx \nonumber\\
&~~~~~ -\int_{\mathcal{D}}\delta\, \big(v^i -\pi_{h_i}(v^i)\big)\,dx \Big|\nonumber\\
&~\leq~ \frac{1}{2}\,\big|\mathcal{A}(v^i; u^i, u^i-u^{i-1})- \mathcal{A}_{h_i}(v^i; u^i, u^i-u^{i-1})\big|+\Big|\int_{\mathcal{D}}\delta\, \big(v^i -\pi_{h_i}(v^i)\big)\,dx\Big|\nonumber\\
&~~~~~+\Big|\int_{\mathcal{D}}\frac{1}{2}\,\Big(\dfrac{( 1- \kappa)\,\nabla u^i\cdot \nabla u^{i-1} }{\big[ 1 + \beta^{\alpha}\,\big((1-\kappa)(v^i)^2+\kappa \big)^\alpha\, | \nabla u^i|^{2\alpha}\big]^{\frac{1}{\alpha}+1}}\Big)\times \big((v^i)^2 -\pi_{h_i}(v^i)^2 \big) \, dx \Big|\nonumber\\
&~\leq~ \frac{1}{2}\,\big|\mathcal{A}(v^i; u^i, u^i-u^{i-1})- \mathcal{A}_{h_i}(v^i; u^i, u^i-u^{i-1})\big|+\Big|\int_{\mathcal{D}}\delta\, \big(v^i -\pi_{h_i}(v^i)\big)\,dx\Big|\nonumber\\
&~~~~~+\frac{1}{4}\Big|\int_{\mathcal{D}}\,\Big(\dfrac{( 1- \kappa)\,|\nabla u^i|^2}{\big[ 1 + \beta^{\alpha}\,\big((1-\kappa)(v^i)^2+\kappa \big)^\alpha\, | \nabla u^i|^{2\alpha}\big]^{\frac{1}{\alpha}+1}} \Big)\times \big((v^i)^2 -\pi_{h_i}(v^i)^2 \big) \, dx \Big| \nonumber\\
&~~~~~+\frac{1}{4}\Big|\int_{\mathcal{D}}\,\Big(\dfrac{( 1- \kappa)\,|\nabla u^{i-1}|^2}{\big[ 1 + \beta^{\alpha}\,\big((1-\kappa)(v^i)^2+\kappa \big)^\alpha\, | \nabla u^i|^{2\alpha}\big]^{\frac{1}{\alpha}+1}} \Big)\times \big((v^i)^2 -\pi_{h_i}(v^i)^2 \big) \, dx \Big|\nonumber\\
&~\leq~ \frac{1}{2}\,\big|\mathcal{A}(v^i; u^i, u^i-u^{i-1})- \mathcal{A}_{h_i}(v^i; u^i, u^i-u^{i-1})\big|+\frac{1}{4}|\mathcal{J}{(u^i, v^i)}-\mathcal{J}_{h_i}{(u^i, v^i)}|\nonumber\\
&~~~~~+\frac{1}{4}\Big|\int_{\mathcal{D}}\,\Big(\dfrac{( 1- \kappa)\,|\nabla u^{i-1}|^2}{\big[ 1 + \beta^{\alpha}\,\big((1-\kappa)(v^i)^2+\kappa \big)^\alpha\, | \nabla u^i|^{2\alpha}\big]^{\frac{1}{\alpha}+1}}  \Big)\times \big((v^i)^2 -\pi_{h_i}(v^i)^2 \big) \, dx \Big| \nonumber \\
&~~~~~+\Big|\int_{\mathcal{D}}\delta\, \big(v^i -\pi_{h_i}(v^i)\big)\,dx\Big|.
\end{align}
Thus, we conclude that
\begin{align}
&|\mathcal{J}{(u^i, v^i)}-\mathcal{J}_{h_i}{(u^i, v^i)}|\leq \frac{2}{3}\,\big|\mathcal{A}(v^i; u^i, u^i-u^{i-1})- \mathcal{A}_{h_i}(v^i; u^i, u^i-u^{i-1})\big|\nonumber\\
&~~~~+\frac{1}{3}\,\Big|\int_{\mathcal{D}}\,\Big(\dfrac{( 1- \kappa)\,|\nabla u^{i-1}|^2\,\big((v^i)^2 -\pi_{h_{i-1/2}}(v^i)^2 \big)}{\big[ 1 + \beta^{\alpha}\,\big((1-\kappa)(v^i)^2+\kappa \big)^\alpha\, | \nabla u^i|^{2\alpha}\big]^{\frac{1}{\alpha}+1}} \Big) \, dx \Big|+\frac{4}{3}\Big|\int_{\mathcal{D}}\delta\, \big(v^i -\pi_{h_i}(v^i)\big)\,dx\Big|\nonumber\\
&:= I_5+I_6+I_7. \label{4.42jhinq}
\end{align}
Now we estimate the terms $I_5$, $I_6$, and $I_7$ separately. To estimate $I_5$, we apply the same estimation approach used for $I_2$ in \eqref{4.18I2} to the equation \eqref{4.16Abd}. Through the substitution of $\psi =u^i-u^{i-1}$ into the residual estimate $\psi \mapsto \mathcal{A}(v^i; u^i, \psi)$, we obtain
\begin{align}
I_5&\leq~c_{51}\,\widetilde{\widetilde{\eta}}_{h_i}(u^i,v^i)\times \|\nabla (u^i-u^{i-1})\|_{L^2(\mathcal{D})},
\end{align}
where $c_{51}=2\,c_{45}/3$ and 
\begin{align*}
\widetilde{\widetilde{\eta}}_{h_i}(u^i,v^i)=~\Big\{ \sum_{ \tau_i\in \mathscr{T}_h}\|\nabla v^i\|^4_{L^\infty(\tau_i)}\,\int_{\tau_i} h_{\tau_i}^4\,\Big|\dfrac{(1-\kappa)\, \nabla u^i }{\big[ 1 + \beta^{\alpha}\,\big((1-\kappa)(v^i)^2+\kappa \big)^\alpha | \nabla u^i |^{2\alpha} \big]^{\frac{1}{\alpha}+1}} \Big|^2\,dx \Big\}^{\frac{1}{2}}. \nonumber
\end{align*}
Following the same methodology as in property (iv), we obtain the estimate for $I_6$, as
\begin{align*}
I_6 &\leq \big|\mathcal{B}(u^{i-1}; v^i, v^i)- \mathcal{B}_{h_{i-1/2}}(u^{i-1}; v^i, v^i)\big|\nonumber \\
&\leq c_{52}\, \widehat{\widehat{\eta}}_{h_{i-1/2}}(u^i,v^i) \times \|\nabla v^i\|, 
\end{align*}
where $c_{52}=\frac{1}{3}\,\max\{c_\varrho, c_{eq}, c_{45}\}$, and $\widehat{\widehat{\eta}}_{h_{i-1/2}}(u^i,v^i)$ can be defined similar to equation \eqref{etahathat}. The last term $I_7$ is estimated using the Lemma \ref{4.3lemmaaprLinfty} with $m=1$, as 
\begin{equation}
I_7\leq c_{53}\,\bar{\bar{\eta}}_{h_i},\nonumber
\end{equation}
where $c_{53}=\frac{4\,\delta}{3}\,c_{45}$ and $\bar{\bar{\eta}}_{h_i}(v^i)=\sum_{\tau_i\in \mathscr{T}_h} h_{\tau_i}\,\|\nabla v^i \|_{L^\infty(\tau_i)}$.
Utilizing the bounds of $I_5$, $I_6$, and $I_7$ in \eqref{4.42jhinq}, we arrive at 
\begin{align*}
|\mathcal{J}{(u^i, v^i)}-\mathcal{J}_{h_i}{(u^i, v^i)}|&\leq c_{51}\,\widetilde{\widetilde{\eta}}_{h_i}(u^i,v^i)\times \|\nabla (u^i-u^{i-1})\|_{L^2(\mathcal{D})}+ c_{52}\, \widehat{\widehat{\eta}}_{h_{i-1/2}}(u^i,v^i) \times \|\nabla v^i\| \nonumber\\
&~~~+c_{53}\,\bar{\bar{\eta}}_{h_i}(v^i).
\end{align*}
Note that $\lim_{i \rightarrow \infty} \Xi_{(RF)_i}=0$ implies the convergence of $\widetilde{\widetilde{\eta}}_{h_i},\, \widehat{\widehat{\eta}}_{h_{i-1/2}}$ and $\bar{\bar{\eta}}_{h_i}$ to zero as $i \rightarrow \infty$. Hence, we have
\begin{align*}
|\mathcal{J}{(u^i, v^i)}-\mathcal{J}_{h_i}{(u^i, v^i)}|\longrightarrow~ 0 \quad \text{as}  \quad i \rightarrow \infty.
\end{align*}
Thus, the property $(v)$ is complete. Finally, we establish the last property. \smallskip

\noindent
\textbf{(vi) \tt Energy descent condition.} Invoking property (iii) assures that $\mathcal{J}_{h_i}{(u^i, v^i)} \leq~ \mathcal{J}_{h_{i-1/2}}{(u^i, v^i)}$, leading to
\begin{align*}
\mathcal{J}(u^i,v^i)&\leq~ \mathcal{J}_{h_i}{(u^i, v^i)}+|\mathcal{J}{(u^i, v^i)}-\mathcal{J}_{h_i}{(u^i, v^i)}|\nonumber\\
&\leq~\mathcal{J}{(u^{i-1}, v^i)}+|\mathcal{J}_{h_{i-1/2}}{(u^i, v^i)}-\mathcal{J}{(u^{i-1}, v^i)}|+|\mathcal{J}{(u^i, v^i)}-\mathcal{J}_{h_i}{(u^i, v^i)}|.
\end{align*}
Applying the properties $(iv)$ and $(v)$, this produces
\begin{align*}
\liminf_{i\rightarrow \infty } \mathcal{J}(u^i,v^i)\leq~ \liminf_{i\rightarrow \infty } \mathcal{J}(u^{i-1},v^i)
\end{align*}
An analogous argument for the right-hand inequality finalizes the proof.
\end{proof}
The sequence $\{(u_k, v_k)\}_{k=1}^\infty \in \mathbb{V}_f\times \mathbb{V}_c$ produced by {\tt Algorithm-2} in section \ref{adaptalgo2} yields the following theorem, which proves its convergence under the following hypothesis:
``{\tt Hypothesis (H):} Assuming that {\tt Step 2} in {\tt Algorithm-2} terminates in a finite number of iterations."
With no termination criterion, the sequence $\{(u_k, v_k)\}_{k=1}^\infty$ converges to a critical point $(u,v)$ in $\mathbb{V}_f\times \mathbb{V}_c^\infty$.
\begin{theorem}[Convergence without Termination] \label{Cngwithouttol}
Let $\mathcal{D}\subset \mathbb{R}^d$ be an open bounded domain. Further, we assume that there exists a sequence $\{(u_k,v_k)\}_{k=1}^{\infty}$ in $\mathbb{V}_f\times \mathbb{V}_c$ with $v_k(x)\in [0, 1]$ for $a.e.\; x\in \mathcal{D}$, generated by the \textbf{\tt Modified Adaptive Algorithm-3} 
 (Sec \ref{madpiii}) under the {\tt Hypothesis \textbf{\tt (H)}}, and also satisfied the condition \eqref{4.30abinq}, Theorem \ref{Cngwithtol}.
Then, there exists a subsequence $\{(u_{k_n}, v_{k_n})\}_{n=1}^{\infty}$ of $\{(u_k, v_k)\}_{n=1}^{\infty}$ and a pair $(u,v)$ in $\mathbb{V}_f\times \mathbb{V}_c$ with $v(x)\in [0, 1]$  $a.e.\; x\in \mathcal{D}$, such that $u_{k_n}$ and $v_{k_n}$ converges strongly to $u$ and $v$ as $n \rightarrow \infty$, respectively, in $\mathbb{V}$. Additionally, $u$ and $v$ satisfy 
\begin{align} \label{4.6criticalpts}  
&\mathcal{A}(v; u, \psi) =0  \quad \; \forall \,   \psi\in \mathbb{V}_d,\quad \text{and} \quad 
\mathcal{B}(u; v, \varphi)=0  \quad \; \forall \,   \varphi \in \mathbb{V}_c^\infty. 
\end{align}
That is, the function $\mathcal{J}(\cdot,\cdot)$ has a critical point $(u,v)$ in $\mathbb{V}_f\times \mathbb{V}_c^\infty$.
\end{theorem}
\begin{proof}
The proof of the theorem executed by the following steps. In the first step, we demonstrate the existence of a convergent subsequence
$\{(u_{i_n},v_{i_n})\}_{n=1}^{\infty}$ of $\{(u_i,v_i)\}_{i=1}^{\infty}$. \\

\noindent
\textbf{\tt Step 1.} According to Step 1 of Theorem \ref{Cngwithtol}, the sequence $\{(u_i,v_i)\}_{i=1}^{\infty}$ is bounded in $\mathbb{V}_f \times \mathbb{V}_c$. As a result, it permits a weakly convergent subsequence $\{(u_{i_n},v_{i_n})\}_{n=1}^{\infty}$, which converges to $(u,v)$ as $n \rightarrow \infty $ in $\mathbb{V} \times \mathbb{V}$. Thus, we acquire $(u_{i_n},v_{i_n}) \overset{w}{\longrightarrow} (u, v),\, n \rightarrow \infty$, where $u\in \mathbb{V}_f$, $v\in \mathbb{V}_c$, and $0\leq v(x) \leq 1$ a.e. in $\mathcal{D}$. However, the careful labeling of this subsequence is crucial for the following proof. 

Next, we shall demonstrate that the expression $\mathcal{A}(v; u, \psi)$ vanishes for every $\psi\in\mathbb{V}_d$.\smallskip

\noindent
\textbf{\tt Step 2.} According to Theorem \ref{Cngwithtol}, the sequence $\{(u_{i_n},v_{i_n})\}_{n=1}^{\infty}$ accomplishes the following
\begin{align*}
&\mathcal{A}(v_{i_n}; u_{i_n}, \psi)\leq~ \gamma^1_{i_n}\|\nabla \psi\|_{L^2(\mathcal{D})}, \quad \forall  n\in \mathbb{N} \quad \text{and } \forall \psi \in \mathbb{V}_d 
\end{align*}
with $\gamma^1_{i_n}\rightarrow 0$ as $n \rightarrow \infty$. It is easy to demonstrate that, using the same logic as in the proof of \textbf{\tt Step 2} of Theorem \ref{Cngwithtol},
$$\mathcal{A}(v; u, \psi) =~0, \; \forall \,   \psi\in \mathbb{V}_d.$$
Furthermore, the subsequence $\{\nabla u_{i_n}\}_{n=1}^\infty$ converges strongly to $\nabla u$ as $n \rightarrow \infty$ in $(L^2(\mathcal{D}))^d$, that is, $\nabla u_{i_n} \longrightarrow \nabla u$ as $n \rightarrow \infty$. 

Consider $\{v_{i_n+1}\}_{n=1}^\infty$ as a subsequence of $\{v_i\}_{i=1}^\infty$ that converges weakly to some $\bar{v}$ in $\mathbb{V}$. In the following step, we will show that this subsequence strongly converges to $\bar{v}$ in $\mathbb{V}$ as $n \rightarrow \infty$.\smallskip

\noindent
\textbf{\tt Step 3.} We know that $\mathbb{V}$ is compactly embedded in $L^2(\mathcal{D})$, hence, the subsequence $\{v_{i_n+1}\}_{n=1}^\infty$ strongly converges to $\bar{v}$ in $L^2(\mathcal{D})$. We want to demonstrate that as $n\rightarrow \infty$, $\|\nabla v_{i_n+1}\|_{L^2(\mathcal{D})} \longrightarrow \|\nabla \bar{v}\|_{L^2(\mathcal{D})}$. This, along with weak convergence, will lead to strong convergence.

Additionally, it should be noted that the sequence $\{(u_{i_n}, v_{i_n+1})\}_{n=1}^\infty$ complies with the inequality \eqref{4.30abinq}, Theorem \ref{Cngwithtol}, thus
\begin{align*}
&\mathcal{B}(u_{i_n}; v_{i_n+1}, \varphi)\leq~ \gamma^2_{i_n}\,\|\nabla \varphi \|_{L^2(\mathcal{D})}, \quad \forall n\in \mathbb{N} \quad \text{ and } \; \forall \varphi \in \mathbb{V}_c. 
\end{align*}
By substituting $\varphi= v_{i_n+1}-\bar{v} \in \mathbb{V}_c$ into the bilinear form $\mathcal{B}(u_{i_n}; v_{i_n+1}, \varphi)$ and rearranging terms, we have
\begin{align}
&2\,\rho\,\int_{\mathcal{D}} |\nabla v_{i_n+1}|^2\, dx  \leq~ \gamma^2_{i_n}\,\|\nabla \varphi \|_{L^2(\mathcal{D})}+2\,\rho\,\int_{\mathcal{D}} \nabla v_{i_n+1} \cdot \nabla \bar{v} \, dx \nonumber\\
&~~~+ \int_{\mathcal{D}}\delta\,(v_{i_n+1}-\bar{v})\, dx + \int_{\mathcal{D}}\dfrac{(1-\kappa)\,[|\nabla u|^2-|\nabla u_{i_n}|^2]\;v_{i_n+1} \,(v_{i_n+1}-\bar{v})}{\left( 1 + \beta^{\alpha}\,\big((1-\kappa)\, v_{i_n+1}^2+\kappa\big)^\alpha |\nabla u_{i_n}|^{2\alpha}\right)^{\frac{1}{\alpha}+1}} \, dx\nonumber\\
&~~~-\int_{\mathcal{D}}\dfrac{(1-\kappa)\,|\nabla u|^2\;v_{i_n+1} \,(v_{i_n+1}-\bar{v})}{\left( 1 + \beta^{\alpha}\,\big((1-\kappa)\, v_{i_n+1}^2+\kappa\big)^\alpha |\nabla u_{i_n}|^{2\alpha}\right)^{\frac{1}{\alpha}+1}} \, dx.\nonumber
\end{align}
Utilizing $\bar{v},\, v_{i_n+1}\in [0, 1]$, to obtain
\begin{align}
&2\,\rho\,\int_{\mathcal{D}} |\nabla v_{i_n+1}|^2\, dx  \leq~ \gamma^2_{i_n}\,\|\nabla \varphi \|_{L^2(\mathcal{D})}+2\,\rho\,\Big|\int_{\mathcal{D}} \nabla v_{i_n+1} \cdot \nabla \bar{v} \, dx \Big| + \delta\,\int_{\mathcal{D}} |v_{i_n+1}-\bar{v}|\, dx \nonumber\\
&~~~+ |1-\kappa|\int_{\mathcal{D}}|\nabla u-\nabla u_{i_n}|\,(|\nabla u|+|\nabla u_{i_n}|)\, dx+|1-\kappa|\,\int_{\mathcal{D}}|\nabla u|^2\,|v_{i_n+1}-\bar{v}|\, dx.\nonumber
\end{align}
Since $|v_{i_n+1}|,\, 1/\left( 1 + \beta^{\alpha}\,\big((1-\kappa)\, v_{i_n+1}^2+\kappa\big)^\alpha \|\nabla u_{i_n}\|^{2\alpha}\right)^{\frac{1}{\alpha}+1}  \leq 1,$ and $\gamma^2_i \longrightarrow 0$ as $n \rightarrow \infty$, this leads to
\begin{align}
2\,\rho\,\limsup_{n\rightarrow \infty} \int_{\mathcal{D}} |\nabla v_{i_n+1}|^2\, dx  \leq ~ 2\,\rho\, \int_{\mathcal{D}} | \nabla \bar{v}|^2 \,dx. 
\end{align}
Taking use of the  weak lower semi-continuity, we infer
\begin{align}
\int_{\mathcal{D}} | \nabla \bar{v}|^2 \,dx \leq~ \liminf_{n\rightarrow \infty} \int_{\mathcal{D}} |\nabla v_{i_n+1}|^2\, dx \leq \limsup_{n\rightarrow \infty} \int_{\mathcal{D}} |\nabla v_{i_n+1}|^2\, dx  \leq ~ \int_{\mathcal{D}} | \nabla \bar{v}|^2 \,dx. 
\end{align}
Thus, the subsequence $\{\|\nabla v_{i_n+1}\|_{L^2(\mathcal{D})}\}_{n=1}^\infty$ converges to $\|\nabla \bar{v}\|_{L^2(\mathcal{D})}$. Combining this with weak convergence, we achieve that $\{v_{i_n+1}\}_{n=1}^\infty$ strongly converges to $\bar{v}$ in $\mathbb{V}$.  This completes the proof.

We now show that, for all $\varphi \in \mathbb{V}_c^\infty$, $\mathcal{B}(u; v, \varphi)=0$ if $v=\bar{v}$.\smallskip

\noindent
\textbf{\tt Step 4.} Notably, the sequence $\{(u_{i_n}, v_{i_n+1})\}_{n=1}^\infty$ fulfills the following inequality \eqref{4.30abinq}, and hence
\begin{align*}
&\mathcal{B}(u_{i_n}; v_{i_n+1}, \varphi)\leq~ \gamma^2_{i_n}\,\|\nabla \varphi \|_{L^2(\mathcal{D})}, \quad \forall n\in \mathbb{N} \quad \text{and} \quad \varphi \in \mathbb{V}_c^\infty. 
\end{align*}
By adopting the similar strategy as \textbf{\tt Step 2-(ii)} of Theorem \ref{Cngwithtol} with $v=\bar{v}$ and letting $n\rightarrow \infty$,  we determine
$$\mathcal{B}(u; \bar{v}, \varphi)=~0, \; \forall\,   \varphi \in \mathbb{V}_c^\infty.$$
It remains to prove that $v=\bar{v}$. Property $(vi)$ of Lemma \ref{410lmadpalg} immediately leads to this equality, and as a result
\begin{align}
\mathcal{J}(u,v)\leq~ \liminf_{n \rightarrow \infty} \mathcal{J}(u_{i_n+1},v_{i_n+1})\leq~ \liminf_{n \rightarrow \infty} \mathcal{J}(u_{i_n},v_{i_n+1})= \mathcal{J}(u,\bar{v}). \label{450inqjc}
\end{align}
Note that $\bar{v}$ is a critical point of the strictly convex functional  $\mathcal{J}(u,v)$. It confirms that $\bar{v}$ is a unique minimizer of $\mathcal{J}$. Consequently, $\mathcal{J}(u,v)\leq \mathcal{J}(u,\bar{v})$ for all $v$, with equality if and only if $v=\bar{v}$. Therefore, we infer that $(u,\bar{v})$ is the only critical point of $\mathcal{J}(\cdot,\cdot)$, which concludes the proof of the theorem.
\end{proof}
\begin{remark}
The physical parameters $\alpha$, $\beta$, $\rho$, and $\delta$ in Theorem \ref{Cngwithouttol} are chosen to ensure the strict convexity of the functional $\mathcal{J}(u,v)$. This selection is crucial for guaranteeing the uniqueness of the minimizer; otherwise, Step 4 of the proof may not hold.
\end{remark}
\section{Numerical Assessments}
This section explores the theoretical implications of our work, illustrated by a concrete example employing {adaptive algorithms}. As detailed in the previous sections, this example focuses on an {elastic unit-square domain} featuring a {single edge-crack} subjected to {antiplane shear boundary loading}. The computational aspects of this study were meticulously implemented using an {adaptive finite element code}, custom-developed by the authors in {\textsf{C++}} using an open-source {\textsf{deal.II} library} \cite{arndt2021deal}.

For this specific example, we calculated the critical parameters, denoted as $u$ and $v$, using two distinct adaptive algorithms: {\tt Adaptive Algorithm-1} and {\tt Adaptive Algorithm-2}. These algorithms are thoroughly defined and discussed in Subsections \ref{algorithm1withtol} and \ref{adaptalgo2}, respectively. The computational domain for this analysis is a {rectangular region} $\mathcal{D}=[0, 1]\times [0, 1]$. A key feature of this domain is a {slit} originating from the point $(0.5, 1)$, as visually represented in Figure \ref{crackdomain}.
\begin{exam}[A crack domain] \label{exm1}
We consider a rectangular domain $\mathcal{D}=[0, 1]\times [0, 1]$ with a slit eminating from $(0.5, 1)$ which shown in the following Figure \ref{crackdomain}. The function $f(x,t)$ represents the incremental anti-plane displacement given by 
\begin{align*}
f(x,t)=
\begin{cases}
-c\,t    \quad  \text{on} \quad (0, 1)\times (0.5, 1), \\
~~ c\, t  \quad  ~\text{on} \quad (0.5, 1)\times (1, 1).
\end{cases}
\end{align*}
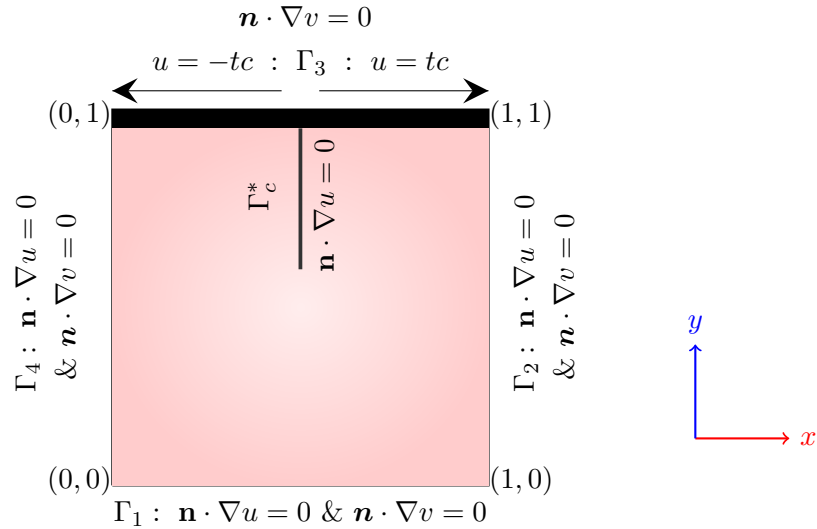
\begin{figure}[H]
\centering
\begin{tikzpicture}[scale=1.25]
\tikzset{myptr/.style={decoration={markings,mark=at position 1 with %
{\arrow[scale=3,>=stealth]{>}}},postaction={decorate}}}
\filldraw[draw=black, thick] (0,0) -- (4,0) -- (4,4) -- (0,4) -- (0,0);
\shade[inner color=red!7, outer color=red!20] (0,0) -- (4,0) -- (4,3.8) -- (0,3.8) -- (0,0);
\node at (-0.35,0.05)   {$(0,0)$};
\node at (4.35,0.05)   {$(1,0)$};
\node at (-0.35,3.95)   {$(0,1)$};
\node at (4.35, 3.95)   {$(1,1)$};
\node at (-0.9, 3.2)[anchor=east, rotate=90]{$\Gamma_4:~ \mathbf{n}\cdot \nabla u=0$};
\node at (-0.5, 3)[anchor=east, rotate=90]{$\&~ \boldsymbol{n} \cdot \nabla v=0$};
\node at (4.4, 3.2)[anchor=east, rotate=90]{$\Gamma_2:~ \mathbf{n}\cdot \nabla u=0$};
\node at (4.8, 3)[anchor=east, rotate=90]{$\&~ \boldsymbol{n} \cdot \nabla v=0$};
\node at (2, -0.3) {$\Gamma_{1}:~ \mathbf{n}\cdot \nabla u=0  ~\&~ \boldsymbol{n} \cdot \nabla v=0$};
\node at (2, 4.5) {$ u=-tc ~: ~\Gamma_{3} ~:~ u=tc$};
\node at (2, 5) {$ ~ \boldsymbol{n} \cdot \nabla v=0$};
\draw [myptr](1.8, 4.2)--(0.0, 4.2);
\draw [myptr](2.2, 4.2)--(4.0, 4.2);
\draw [line width=0.5mm, black!80]  (2,3.8) -- (2,2.3);
\node at (1.59,3.35)[anchor=east, rotate=90]{$\Gamma^*_c$};
\node at (2.25,3.8)[anchor=east, rotate=90]{$\mathbf{n}\cdot \nabla u=0$};
    \def\xOrigin{6.2}
    \def\yOrigin{0.5} 
    \draw[->, color=red, thick] (\xOrigin, \yOrigin) -- (\xOrigin + 1, \yOrigin) node[right] {$x$};
    \draw[->, color=blue, thick] (\xOrigin, \yOrigin) -- (\xOrigin, \yOrigin + 1) node[above] {$y$};
\end{tikzpicture}
\caption{A domain and the boundary indicators.~~~~~~~~~~~~~~~~~~~~~~}\label{crackdomain}
\label{fig:h-conv}
\end{figure}
\end{exam}
For the numerical calculations in Example \ref{exm1}, we simulated 60 time steps, each with a uniform step size of $k=0.01$. The adaptive algorithms employed specific tolerances: $\Xi_{RF}=0.01$, $\Xi_{CR}=10^{-04}$, $\Xi_{v}=10^{-04}$, and $\Xi_{v_n}=10^{-06}$. During the marking strategy for mesh refinement, the parameter $\vartheta$ was set to $0.5$. The following parameters governed the overall computation: $\alpha=1.0$, $\beta=1.0$, $\kappa=10^{-10}$, and $\epsilon=10{h_\tau}$. The irreversibility criterion was established with a parameter of $10^{-2}$. Additionally, $\lambda_c$ and $c_w$ were set to $2.7$ and $8/3$, respectively. In this example, we determined the final field $v$, bulk energy, surface energy, and total energy using Adaptive Algorithms -1  and 2, as detailed in Subsections \ref{algorithm1withtol} and \ref{adaptalgo2}. These algorithms also generated the final computational meshes. The nonlinear problem for the mechanics was solved using Picard's iteration technique, and we deliberately omitted condition \eqref{aij4.2}. At each time step, the initial fracture field $v$ was set to the final $v$ at the preceding time step. The sole exception was the first time step, where the initial fracture field $v$ was initialized to $1.0$.
\begin{figure}[H] 
    \centering
    \begin{subfigure}{0.49\textwidth} 
        \includegraphics[width=1.0\linewidth]{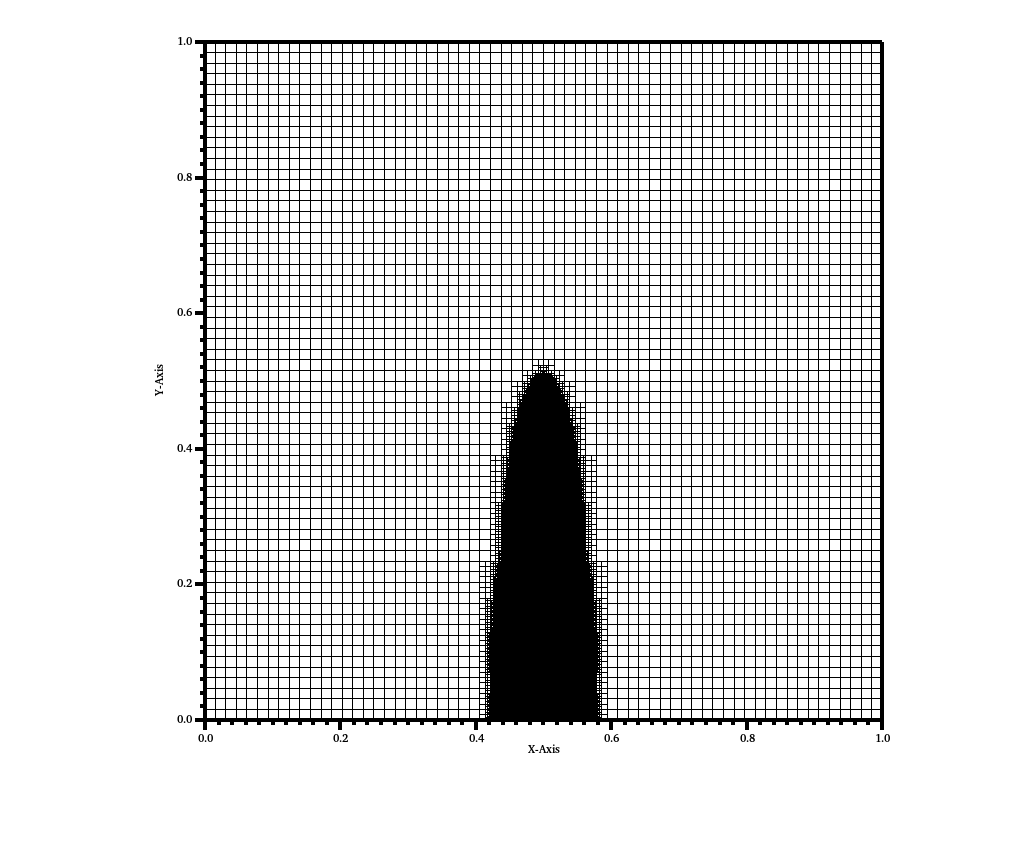} 
        \caption{Final mesh from algorithm-1}
        \label{fig:fig1}
    \end{subfigure}
    \hfill 
    \begin{subfigure}{0.49\textwidth}
        \includegraphics[width=1.0\linewidth]{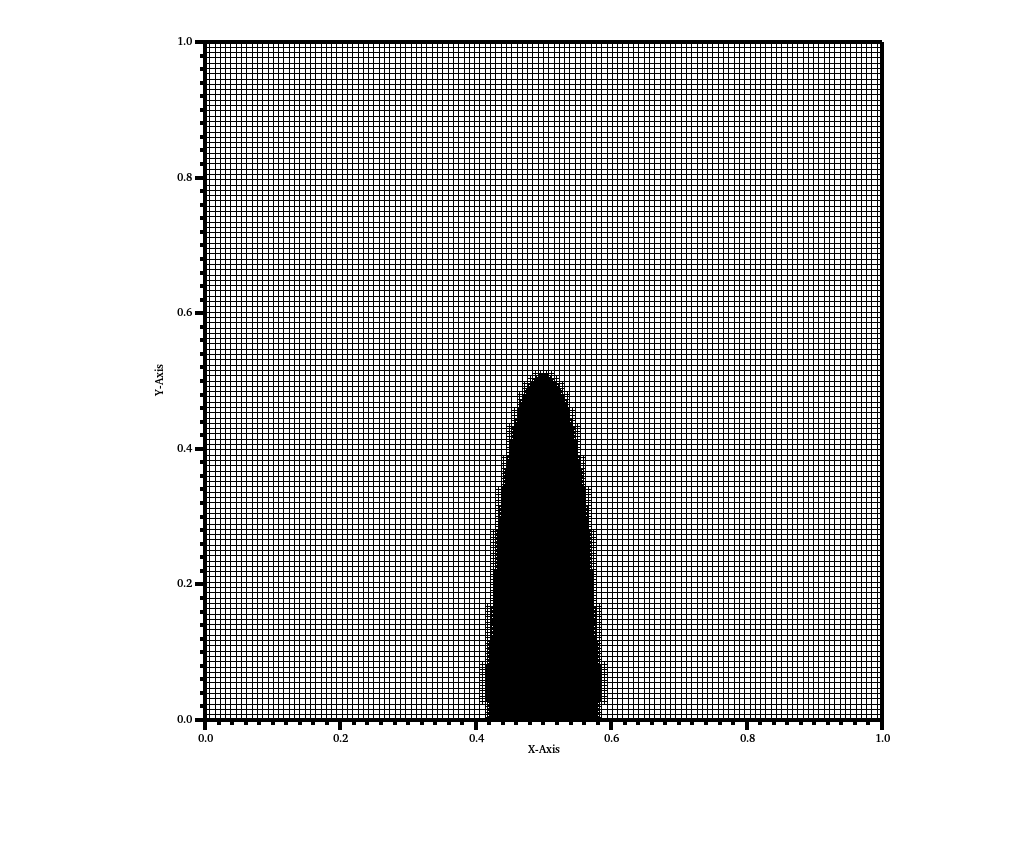} 
        \caption{Final mesh from algorithm-2}
        \label{fig:fig2}
    \end{subfigure}
\vspace{\baselineskip} 
 \centering
    \begin{subfigure}{0.49\textwidth} 
        \includegraphics[width=1.0\linewidth]{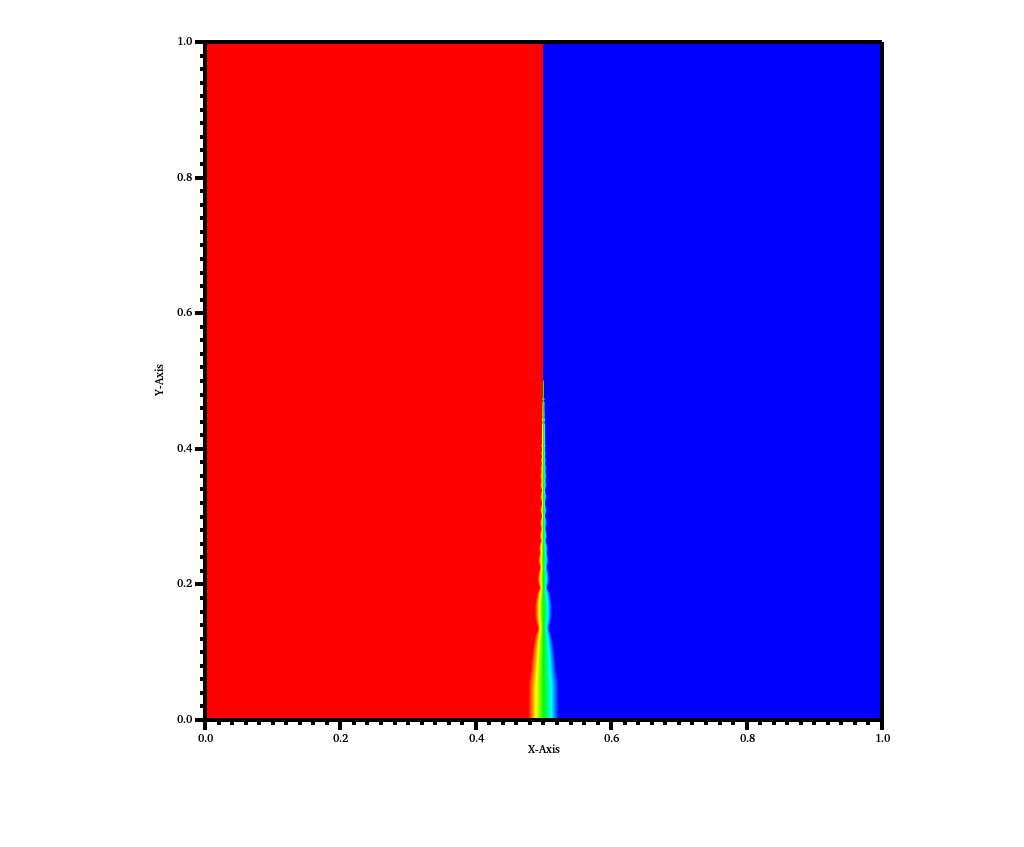} 
        \caption{Final $u$-field from algorithm-1}
        \label{fig:fig1}
    \end{subfigure}
    \hfill 
    \begin{subfigure}{0.49\textwidth}
        \includegraphics[width=1.0\linewidth]{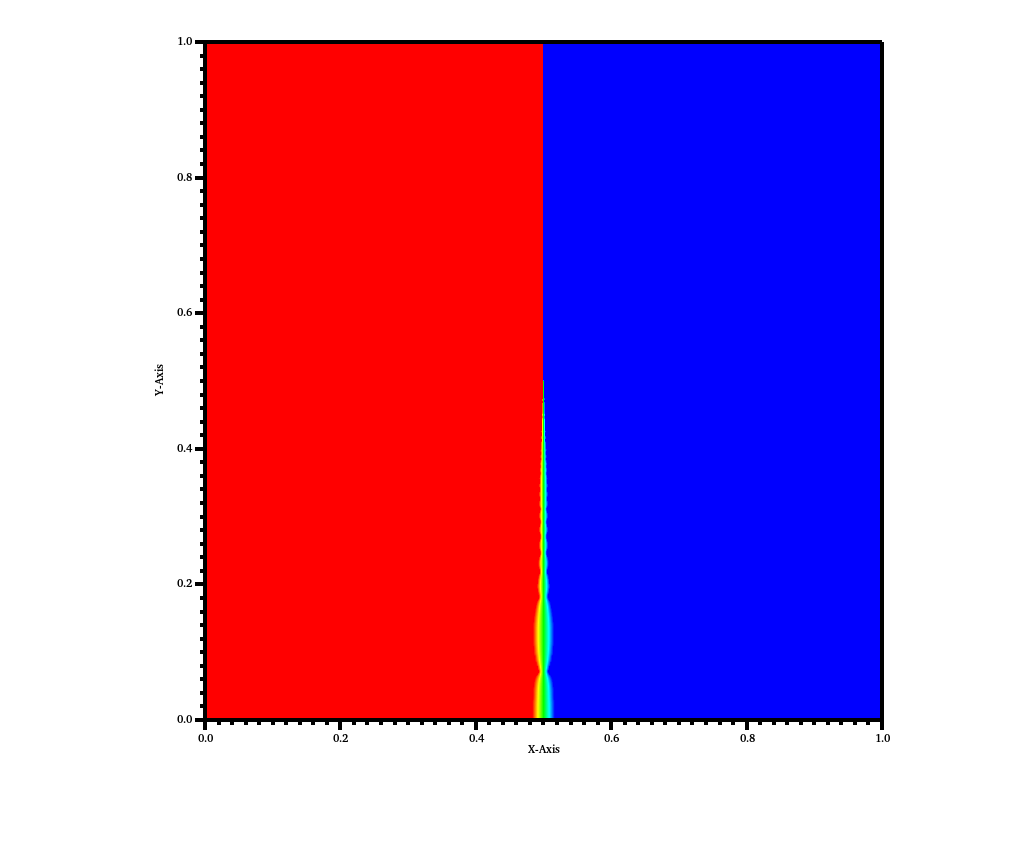} 
        \caption{Final $u$-field from algorithm-2}
        \label{fig:fig2}
    \end{subfigure}
    \caption{Results from Algorithm-1 (left panel) and Algorithm-2 (right panel), illustrating the final computational mesh and final $u$-field. }
    \label{fig1}
 \end{figure} 
The results generated by both Algorithm-1 and Algorithm-2 are visually presented in Figures~\ref{fig1}, and \ref{fig2}.  While a comparable overall solution quality was achieved by both approaches, distinct operational differences were observed. Specifically, Algorithm-1 is characterized by its more frequent alternating minimization steps between the two governing equations. In contrast, Algorithm-2 incorporates a greater number of refinement steps during the iterative solution of the linear system of equations. Consequently, the total number of degrees of freedom at the final computational stage was marginally higher for Algorithm-2 than for Algorithm-1. This outcome is consistent with expectations, as Algorithm-1 processes the $v$-equation on a more refined mesh, which in turn contributes to a better quality solution for the $u$-equation.
\begin{figure}[H]
\vspace{\baselineskip} 
 \centering
      \begin{subfigure}{0.49\textwidth}
      \centering
        \includegraphics[width=1.0\linewidth]{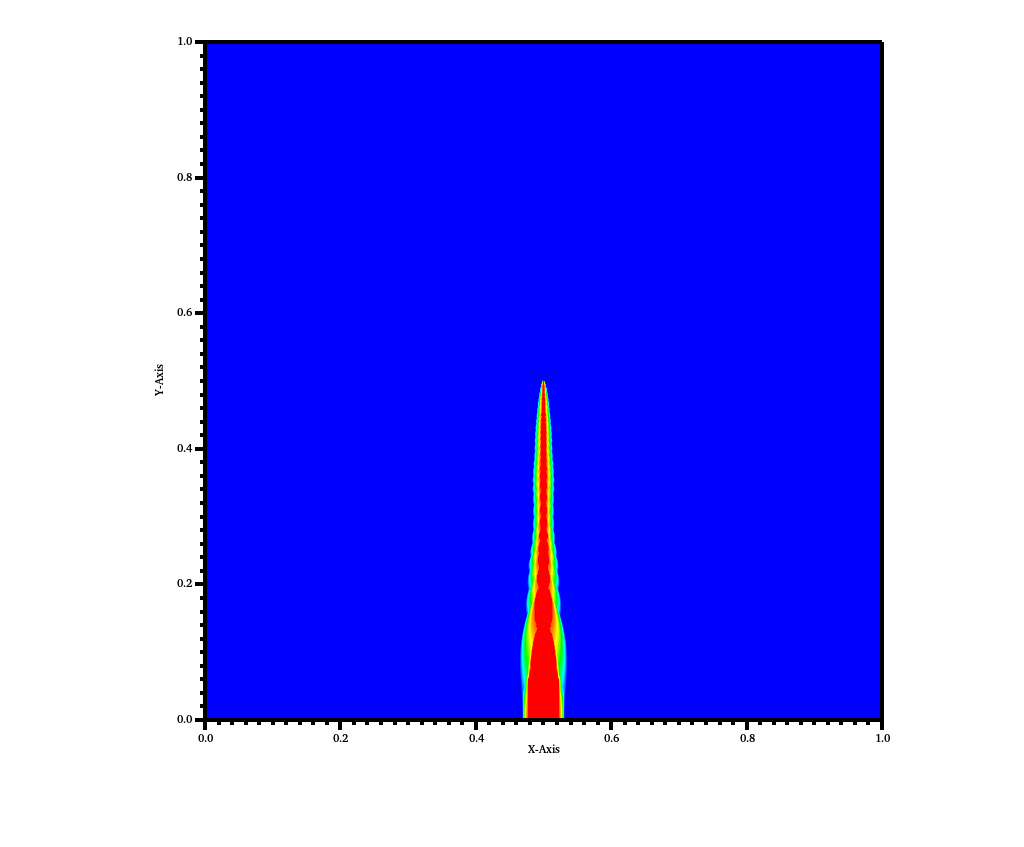} 
        \caption{Final $v$-field from algorithm-1}
        \label{fig:fig3} 
    \end{subfigure}
    \hfill
    \begin{subfigure}{0.49\textwidth} 
        \includegraphics[width=1.0\linewidth]{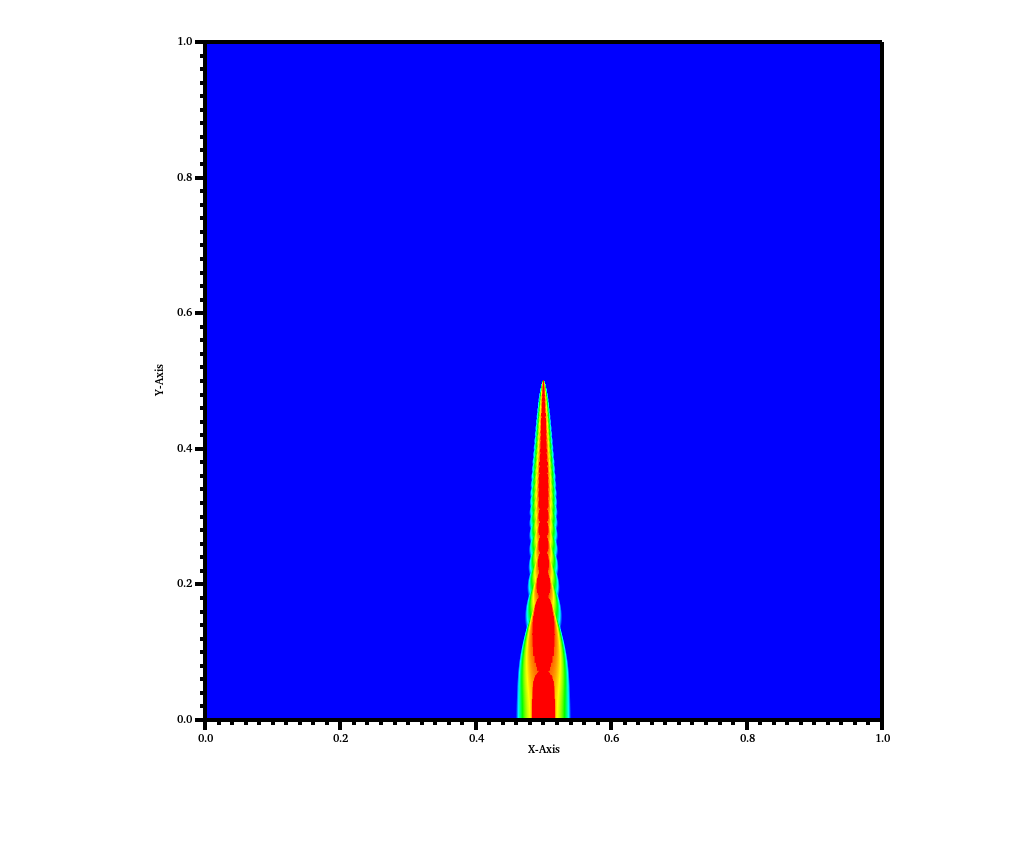} 
        \caption{Final $v$-field from algorithm-2}
        \label{fig:fig4}
    \end{subfigure}
\vspace{\baselineskip} 
 \centering
   \begin{subfigure}{0.49\textwidth}
        \includegraphics[width=0.9\linewidth]{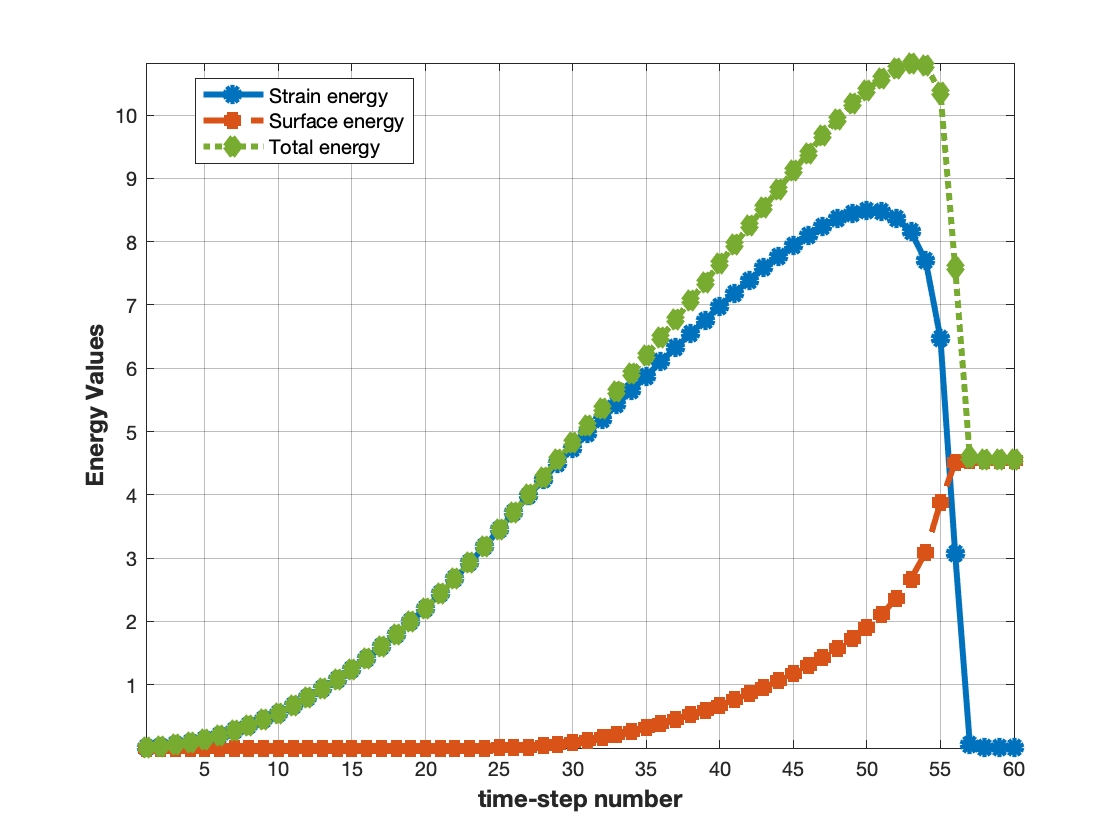} 
        \caption{Energies from algorithm-1}
        \label{fig:fig5}
    \end{subfigure}
    \hfill
    \begin{subfigure}{0.49\textwidth}
        \includegraphics[width=0.9\linewidth]{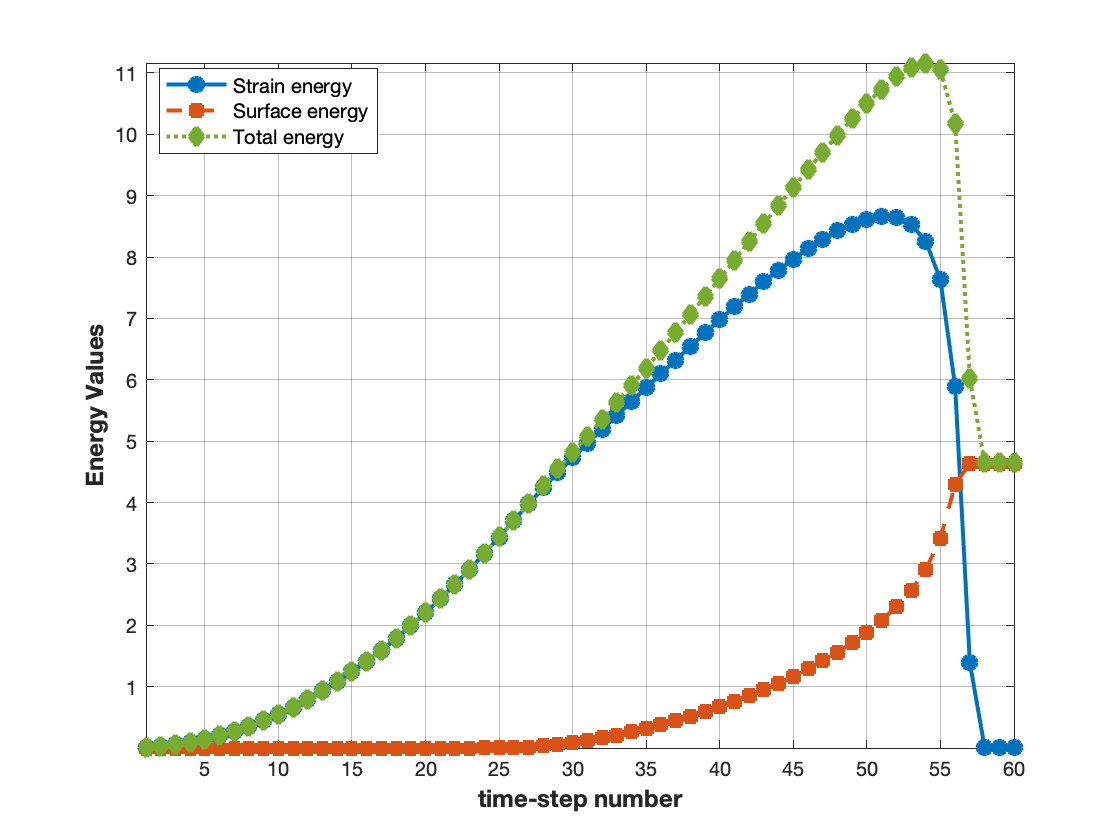} 
        \caption{Energies from algorithm-2}
        \label{fig:fig6}
    \end{subfigure}
    \caption{Results from Algorithm-1 (left panel) and Algorithm-2 (right panel), illustrating the final v-field, and the breakdown of energies (strain, surface, and total). }
    \label{fig2}
\end{figure}
\section{Conclusion} 
In this work, we rigorously analyzed the convergence of adaptive finite element approximations for local minimizers of the Ambrosio-Tortorelli energy functional (AT1), specifically with a nonlinear strain energy density. We centered our investigation on the convergence properties of two distinct adaptive algorithms. A key finding is that these algorithms converge to a critical point of the functional $\mathcal{J}$. We also noted the sensitivity of fracture paths to the chosen algorithmic parameters, emphasizing how parameter selection significantly influences predicted fracture evolution. Additionally, the study underscored the challenges in minimizing the functional $\mathcal{J}_h$ due to its inherent nonlinear and non-convex terms. While identifying global minimizers continues to be an open problem, this research successfully identified local minimizers, providing valuable insights for future studies. These findings advance our understanding of adaptive finite element methods in fracture mechanics. This paves the way for further research into convergence to true minimizers, extending these models to more complex fracture scenarios, and developing even more efficient adaptive algorithms. Our numerical experiments clearly validate the effectiveness of the strategies employed, demonstrating that the sequences generated by the adaptive algorithms consistently converge to a critical point of the functional $\mathcal{J}$, with corresponding residuals simultaneously converging to zero.

This study represents a pivotal initial stride in developing novel adaptive numerical methods for investigating quasi-static crack propagation. Building upon this foundational research, numerous avenues for future work emerge. A similar convergence analysis of finite element algorithms could be extended to address crack-tip evolution in porous elastic solids \cite{yoon2024finite}, orthotropic solids \cite{ghosh2025finite}, and three-dimensional problems \cite{gou2023computational,gou2023MMS,gou2025computational}. Furthermore, the methodologies developed herein could be adapted to analyze three-field formulation problems \cite{fernando2025xi}.
\section{Acknowledgement.}
This material is based upon work supported by the National Science Foundation under Grant No. 2316905. 
\bibliographystyle{plain}
\bibliography{references}
\end{document}